\def \journal{0}

\if\journal1

\RequirePackage{fix-cm}
\documentclass[smallextended]{svjour3}

\else
\documentclass[11pt]{article}

\usepackage{fullpage}
\usepackage{authblk}
\usepackage{amsthm}%
\fi

\usepackage{amsmath,amssymb,amsfonts}
\usepackage[table,xcdraw, dvipsnames]{xcolor}
\usepackage[utf8x]{inputenc}
\DeclareUnicodeCharacter{2212}{-}
\usepackage{lmodern}
\usepackage{rotating}
\usepackage{fancyvrb}
\usepackage{dirtytalk}
\usepackage{cases}
\usepackage{float}
\usepackage{subcaption}
\usepackage{tikz}
\usepackage{hyperref}

\usepackage{graphicx}%
\usepackage{multirow}%
\usepackage{mathrsfs}%
\usepackage[title]{appendix}%
\usepackage{textcomp}%
\usepackage{manyfoot}%
\usepackage{booktabs}%
\usepackage{algorithm}%
\usepackage{algorithmicx}%
\usepackage{algpseudocode}%
\usepackage{listings}%
\usepackage[normalem]{ulem} %

\usepackage{microtype}  %

\usepackage{natbib}
 \bibpunct{[}{]}{,}{n}{}{,}%

\hypersetup{
	colorlinks,
	citecolor=blue,
	filecolor=blue,
	linkcolor=blue,
	urlcolor=black
}

\if\journal0
\newtheorem{theorem}{Theorem}

\newtheorem{example}{Example}%
\newtheorem{property}{Property}%
\newtheorem{lemma}{Lemma}%
\newtheorem{corollary}{Corollary}%
\newtheorem{remark}{Remark}%
\fi

\if\journal1
\smartqed  %
\fi

\newcommand{\korange}{\textcolor[rgb]{0.996, 0.38, 0}}

\def\1{{\bf 1}}
\def\0{{\bf 0}}

\newcommand{\lov}{Lovász}
\newcommand{\sch}{Schrijver}
\newcommand{\R}[1]{\mathbb{R}^{#1}}

\newcommand{\virgo}[1]{``#1''}
\newcommand{\cbra}[1]{\left \{ #1 \right\} }

\newcommand{\ON}{\operatorname}

\newcommand{\M}[1]{%
    \ifx&#1&%
        {\ifmmode M_+(\cdot)\else $M_+(\cdot)$\fi}%
    \else
        {\ifmmode M_+(#1)\else $M_+(#1)$\fi}%
    \fi
}

\newcommand{\N}[1]{%
    \ifx&#1&%
        {\ifmmode N_+(\cdot)\else $N_+(\cdot)$\fi}%
    \else
        {\ifmmode N_+(#1)\else $N_+(#1)$\fi}%
    \fi
}

\newcommand{\Mnod}[2]{\M{\ON{NOD}(#1, #2)}}
\newcommand{\Nnod}[2]{\N{\ON{NOD}(#1, #2)}}
\newcommand{\Mfrac}[1]{\M{\ON{FRAC}(#1)}}
\newcommand{\Nfrac}[1]{\N{\ON{FRAC}(#1)}}
\newcommand{\Mqstab}[1]{\M{\ON{QSTAB}(#1)}}
\newcommand{\Nqstab}[1]{\N{\ON{QSTAB}(#1)}}
\newcommand{\MqstabC}[1]{\M{\ON{QSTAB}(#1, \mathcal{C})}}
\newcommand{\NqstabC}[1]{\N{\ON{QSTAB}(#1, \mathcal{C})}}

\newcommand{\nod}[2]{{\ifmmode \ON{NOD}(#1, #2) \else $\ON{NOD}(#1, #2)$\fi}}
\newcommand{\qstabc}[1]{{\ifmmode \ON{QSTAB}(#1, \mathcal{C}) \else $\ON{QSTAB}(#1, \mathcal{C})$\fi}}
\newcommand{\qstab}[1]{{\ifmmode \ON{QSTAB}(#1) \else $\ON{QSTAB}(#1)$\fi}}

\newif\ifshowdiff
\showdifffalse

\ifshowdiff
    \newcommand{\add}[1]{\kblue{#1}}                      %
    \newcommand{\remove}[1]{\korange{\sout{#1}}}     %
    \newcommand{\replace}[2]{\remove{#1} \add{#2}}   
\else
    \newcommand{\add}[1]{#1}
    \newcommand{\remove}[1]{\ifmmode \else\unskip\fi}
    \newcommand{\replace}[2]{\ifmmode#2\else\unskip#2\fi}
\fi

\newcommand{\ABSTRACT}{
The Lov\'asz theta function $\theta(G)$ provides a very good upper bound on the stability number of a graph $G$. It can be computed in polynomial time by solving a semidefinite program (SDP), which also turns out to be fairly tractable in practice. Consequently, $\theta(G)$ achieves a hard-to-beat trade-off between computational effort and strength of the bound. Indeed, several attempts to improve the theta bound are documented, mainly based on playing around the application of the \N{} lifting operator of Lov\'asz and Schrijver to the classical formulation of the maximum stable set problem \add{(SSP)}. Experience shows that solving such SDPs often struggles against practical intractability and requires highly specialized methods. We investigate the application of such an operator to two different linear formulations \add{of the SSP} based on clique and nodal inequalities, respectively. \replace{Fewer inequalities describe these two and}
{These two formulations are described by fewer inequalities than the natural formulation based on edge inequalities,} yet they guarantee that the resulting SDP bound is at least as strong as $\theta(G)$. Our computational experience, including larger graphs than those previously documented, shows that upper bounds stronger than $\theta(G)$ can be accessed by a reasonable additional effort using the clique-based formulation on sparse graphs and the nodal-based one on dense graphs.
}

\begin{document}

\if\journal1

\title{Application of the Lov\'asz-Schrijver Operator to Compact Stable Set Integer Programs}
\titlerunning{Application of the Lov\'asz-Schrijver Operator to Compact SSP}

\author{Federico Battista \and Fabrizio Rossi \and Stefano Smriglio}

\institute{Federico Battista
             \at
              Department of Industrial and Systems Engineering, Lehigh University, 
              Bethlehem, 18015, PA, USA \\
              \email{feb223@lehigh.edu}           %
           \and
           Fabrizio Rossi    \and
           Stefano Smriglio \at
              DISIM, University of L'Aquila, Via Vetoio, Coppito (AQ), 67100, Italy \\
              \email{\{fabrizio.rossi, stefano.smriglio\}@univaq.it}
}

\date{Received: date / Accepted: date}

\else
\title{Application of the Lov\'asz-Schrijver Operator to Compact Stable Set Integer Programs}

\author[1]{Federico Battista\thanks{\texttt{feb223@lehigh.edu}}}
\author[2]{Fabrizio Rossi\thanks{\texttt{fabrizio.rossi@univaq.it}}}
\author[2]{Stefano Smriglio\thanks{\texttt{stefano.smriglio@univaq.it}}}
\affil[1]{Department of Industrial and Systems Engineering, Lehigh University,
           Bethlehem, PA, 18015, USA }
\affil[2]{DISIM, University of L'Aquila, Coppito (AQ), 67100, Italy}
\fi

\maketitle

\begin{abstract}

\ABSTRACT

\if\journal1


\else


\fi
\end{abstract}

\section{Introduction}
\label{sec:introduction}
Given a simple undirected graph $G=(V,E)$, with vertex set $V$ and edge set $E$, a {\em stable} (or {\em independent}) {\em set} in $G$ is a set of pairwise non-adjacent vertices. A fundamental problem in combinatorial optimization is the \emph{Stable Set Problem}~(SSP), which consists of computing a stable set of maximum cardinality~\cite{GrLoSc2012,NeWo1999}. This value, denoted by $\alpha(G)$, is the {\em stability number} of $G$.
\remove{A weighted version of this problem can be considered if a weight vector $w \in \mathbb{Q}_+^{\replace{n}{|V|}}$ 
is given.}
\add{For the sake of presentation, we do not explicitly address the weighted version of this problem, where a weight vector $w \in \mathbb{Q}_+^{\replace{n}{|V|}}$ is given, although all the results presented remain valid in that setting.}
\add{The} SSP is NP-hard in the strong sense and even hard to approximate~\cite{Ha1996}, it is 
equivalent to the \emph{Maximum Clique} and \emph{Set Packing} problems~\cite{Pa1973} and has 
a wide range of applications in \add{different} practical 
contexts\add{, see, e.g., those illustrated in~\citep{CaSeWi2014, CaGo1999, RoSm2001b} and those surveyed in~\citep{BoBuPaPe1999, WuHa2015}}. \add{The} SSP is naturally formulated as a binary program by the so-called {\em edge inequalities}
\begin{eqnarray}
\alpha\remove{_w}(G)=\max		& \sum_{i \in V} \remove{w_i} x_i & \nonumber \\
\mbox{s.t.}	& \; \; x_i + x_j \le 1, \;	& \{i,j\} \in E, \label{eq:edge} \\
& x_i \in \{0,1\},			& i \in V. \nonumber
\end{eqnarray}

\noindent
As is customary in the literature, we denote by 
$$\ON{FRAC}(G) = \{x \in [0,1]^{|V|}: (\ref{eq:edge}) \mbox{ hold}\}$$ 
the \emph{fractional stable set polytope}, while 
$$\ON{STAB}(G) =\ON{conv} \{x \in \{0,1\}^{|V|} : (\ref{eq:edge}) \mbox{ hold}\}$$ 
denotes the \emph{stable set polytope}, that is, the convex hull of the incidence vectors of all stable sets in $G$.
The intense study of STAB($G$) has led, since the early seventies~\cite{BaPa1976}, to the identification of several classes of valid inequalities, often associated with specific \virgo{facet-producing} sub-graphs~\cite{Bo1998, CaLaMa2000, GrLoSc2012}. In practice, however, linear relaxations based on these inequalities often yielded rather weak upper bounds on the stability number of unstructured graphs. As a consequence, basic branch-and-cut algorithms are not particularly effective~\cite{NeSi1992} and much more sophisticated cut generation procedures, \add{see}~\cite{CoGu2022, CoDeKoMa2018, GiRoSm2013, GiLeRoSm2015, ReOsThRePa2011, RoSm2001}, are required to reduce the performance gap with fast combinatorial exact algorithms, such as those devised in~\cite{Os2002, SaRoJi2011,  ToKa2007}.

\noindent
In his seminal paper, \citet{Lo1979} introduced the celebrated \textit{theta number} of a graph, denoted by $\theta(G)$. It represents an upper bound \replace{for}{on} $\alpha(G)$, which can be computed in polynomial time (up to arbitrary precision) by solving a Semidefinite Program (SDP)\remove{~\cite{GrLoSc2012}}. Since its introduction, $\theta(G)$, along with the associated non-polyhedral convex {\em theta body} $\ON{TH}(G)$\remove{~\cite{GrLoSc2012}}, turned out to be powerful tools. For instance, they allowed to prove that the SSP can be solved in polynomial time when $G$ is a perfect graph\remove{~\cite{GrLoSc2012}}. \add{For the description of these results we refer the reader to the book of~\citet{GrLoSc2012}.} Indeed, the practical relevance of $\theta(G)$ is no less. On one side, it often turns out to be a significantly stronger upper bound than those from linear relaxations, \add{as documented in}~\cite{GiLeRoSm2009, LeMaRoSm2016}. On the other side, the structure of the associated SDP relaxation often mitigates numerical difficulties that typically \replace{affect}{affects} SDP algorithms. As a consequence, $\theta(G)$ realizes a very good trade-off between the quality of the upper bound and computational effort, which makes it attractive to \replace{be incorporated}{incorporate it} within branch-and-bound algorithms, as studied by, e.g., \citet{Wi2009} and more recently by~\citet{GaSiWi2022}.

\noindent
In a successive landmark paper~\citet{LoSc1991} introduced the \N{} {\em lift-and-project} operator, which can be applied to the LP relaxation of any 0-1 LP and returns a stronger SDP relaxation. Its application to $\ON{FRAC}(G)$ has also been intensively investigated and has resulted \add{in a} tighter \add{relaxation} than the Lovász theta relaxation. From a theoretical perspective, its strength has been documented by showing implications of (additional) facet defining inequalities of the stable set polytope (\citet{GiLe2006, LoSc1991}). Moreover,~\citet{BiEsNaTu2017, BiEsNaWa2023} showed that it provides an effective tool to identify other classes of graphs for which the SSP can be solved in polynomial time. Recently, graph classes yielding stable set polytopes with high rank with respect to \N{} have been studied in \cite{YuTunc2024}.

From the computational point of view, things are somewhat more involved. Experiments conducted with \Nfrac{G} \add{by}~\citet{Da2001, BuVa2006} and \add{with} related variants \add{by}~\citet{DuRe2007, GrRe2003}, show that these relaxations may yield stronger upper bounds than $\theta(G)$ but at the price of considerably increasing the computational effort. In fact, solving these large-scale SDPs requires specialized methods, and large instances can hardly be managed.
 Other methods, based on strengthening SDP relaxations by linear inequalities not related to the \N{} operator, have also been presented by~\citet{GaRe2020} and~\citet{Lo2015}. These achieve strong bounds but still require
\remove{a}
significant computational effort.

Overall, we currently face a dichotomic picture, where even minimal progress from the Lov\'asz $\theta$ bound may be paid with a very high computational burden. This paper aims to present new SDP relaxations that can bridge such a gap. These are derived by applying the \N{} operator to alternative LP relaxations. The first relaxation is obtained by replacing edge inequalities with (clearly stronger) clique inequalities. This is unconventional since the linear description may have exponentially many inequalities, making the approach somehow unattractive from theoretical and practical perspectives. However, we show that a natural selection of clique inequalities allows for more compact, tractable, yet stronger, relaxations. An opposite rationale derives the other SDP relaxations: accepting to start from more compact but weaker LP relaxations and relying on the ability of the \N{} operator to recover their potential weakness. Based on this rationale, we investigate a hierarchy of SDP relaxations where the primary level consists of applying the operator to a surrogate relaxation of $\ON{FRAC}(G)$, introduced in~\citet{DeTa1994}, with $O(|V|)$ (so-called \emph{nodal}) inequalities. Coefficient strengthening procedures then derive two progressively stronger relaxations, the first of which takes polynomial time while the second requires the solution of NP-hard (sub-)problems to be built.

 In this work, we first prove that all the proposed SDP relaxations are at least as strong as the Lov\'asz theta relaxation and draw a complete
 theoretical picture by comparing them to each other and to \Nfrac{G}. Then, we discuss implementation issues in handling them and illustrate extensive experiments where the practical strength of relaxations is evaluated and contrasted to the theoretical expectation. We show that significant progress with respect to $\theta(G)$ is achieved in several cases. Relaxations obtained from clique inequalities are more effective for sparse graphs, while those from nodal inequalities are the best as graph density grows. Furthermore, the increase of computational workload required to improve the $\theta$ bound is significantly smoother than that of previous approaches. Indeed, this increased efficiency allowed us experimentation based on significantly larger graphs than previously documented. This experience also sheds new light on the practical potential of the \N{} operator, which emerges more clearly than in the studies documented so far. All the material and code presented is available at~\url{https://github.com/febattista/SDP_lift_and_project} and organized to facilitate the reproduction of relaxations and experiments.

\section{Preliminaries}
 This paper will adopt the following standard notation. The set of integers $\{1 ,...,k\}$ is denoted by $[k]$. The \emph{complement} of a \add{simple} graph $G=(V,E)$ is the graph ${\bar G} = \big( V,{\bar E} \big)$, where ${\bar E} = \big\{ \{i,j\} \subset V: \, \{i,j\} \in (V \times V) \setminus E\add{, i \not= j} \big\}$. Given a vertex $i \in V$, we let $\Gamma_G(i)$ denote the set of \emph{\replace{neighbours}{neighbors}} of $i$ in $G$, i.e., $\Gamma_G(i) = \big\{ j \in V \setminus \{i\}: \{i,j\} \in E \big\}$. When no confusion arises, we let $\Gamma_G(i) = \Gamma(i)$. 
Finally, given a 
\remove{vertex} 
\add{set} $S \subseteq V$, we let $G[S] = \big( S, E(S) \big)$ denote the subgraph of $G$ induced by $S$, and $r(S)$ denote its stability number $\alpha \big( G[S] \big)$ (the so-called \emph{rank} of $S$). We also let ${\cal S}_{n}$ be the set of real symmetric square matrices of order $n$ and ${\cal S}^+_{n}~\subset~{\cal S}_{n}$ denotes the cone of those which are positive semidefinite (PSD). PSDness of a matrix $Y \in {\cal S}^+_{n}$ is also denoted by $Y \succeq 0$. Given any square matrix $X \in \mathbb{R}^{n\times n}$ we denote by \replace{$(X)_{ij}$}{$X_{ij}$} the element at the $i$-th row and the $j$-th column, by $\ON{diag}(X) \in \mathbb{R}^{n}$ its main diagonal and by $\ON{rank}(X)$ its rank. We denote by $\left\langle X,Y\right\rangle = \ON{tr}(XY)$ the standard inner product in ${\cal S}_{n}$, where $\ON{tr}(A)$ is the trace of the square matrix $A$. At last, we denote \replace{with}{by} $J$ the all-ones square matrix and $e$ the all-ones vector; the dimensions should be clear from the context where not explicit.
\add{Given $a \in \R{n}$ and $b \in \R{}$, the linear inequality $a^\top x \leq b$ is \emph{valid} for $\ON{STAB}(G)$ 
if 
$ \left\{ x \in \R{n} \mid a^\top x \leq b \right\} \supseteq \ON{STAB}(G).$}

\subsection{The lift-and-project operator of Lovász and Schrijver}
\label{sec:lift_and_project}
We now review the general statement of the \N{} operator to make the presentation self-consistent, referring the reader to~\cite{CoCoZa2014, Da2001, GoeTunc2001, LoSc1991} for an exhaustive treatment.

Let us consider a 0-1 linear program of the form
\begin{equation}
	\text{max}\left\lbrace c^\top x\ |\ Ax \leq b,\ x_i \in \{0,1\} \text{ for } i \in \replace{I=\{1, \ldots, n\}}{[n]} \right\rbrace, \tag{IP}
	\label{IP}
\end{equation}
where $c \in \mathbb{R}^n$, the inequality system $Ax \le b$ has $A \in \mathbb{R}^{m \times n}$, $b \in \mathbb{R}^m$ and contains the inequalities $0 \le x_i \le 1$, for $i \in [n]$. Each linear inequality of the system $Ax \le b$ is denoted by  $(a^k)^\top x \le b^k$, for $k \in [m]$. The polytope $P=\{x \in \mathbb{R}^n\ |\ Ax \le b\}$ is the linear relaxation of \eqref{IP} and $P_I=\operatorname{conv}\{x \in \{0,1\}^n\ |\ Ax \le b\}$ denotes the convex hull of the feasible solutions to \eqref{IP}. 

\noindent
The application of the {\em lift-and-project} operator to $P$ involves the three following steps:

\subsubsection*{Step 1 (Lifting)} Generate the quadratic inequalities
\begin{align}
	x_i \add{\cdot} ((a^k)^\top x - b^k)  & \leq 0 \quad  i \in [n], k \in [m], \label{constr1} \\
	(1- x_i) \add{\cdot} ((a^k)^\top x - b^k)  & \leq 0 \quad  i \in [n], k \in [m]. \label{constr2}
\end{align}
which are valid for $P$ as $x_i \geq 0$ and $1 - x_i \geq 0$ for all $x \in P$. These inequalities can be linearized through the introduction of the matrix variable $X = xx^\top$ and 
\remove{rewriting}~\eqref{constr1} and~\eqref{constr2} \add{can be rewritten} as\remove{:}

\begin{align}
	\sum_{j \in I} a^k_j X_{ij} - b^kx_i & \leq 0 \quad i \in [n], k \in [m], \label{lin:constr1} \\
	\sum_{j \in I} a^k_j x_{j} - b^k - \sum_{j \in I} a^k_j X_{ij} + b^kx_i & \leq 0 \quad i \in [n], k \in [m], \label{lin:constr2} \\
	X &= xx^\top .\nonumber
\end{align}
\noindent
Then, the quadratic constraint $X = xx^\top$ is relaxed with $X - xx^\top \succeq 0$. By applying the Schur's complement definition $X - xx^\top \succeq 0$ is equivalent to
\begin{equation}
	\begin{pmatrix} 1&x^\top \cr x&X \cr
		\end{pmatrix} \succeq 0. \label{eq:Schur}   
\end{equation}

\subsubsection*{Step 2 (Strengthening)} 

Add to constraints  (\ref{lin:constr1})--(\ref{eq:Schur}) 	

\begin{equation}
	X_{ii} = x_i \qquad i \in [n], \label{lin:constr3}
\end{equation}

\noindent
as any binary vector $x \in P_I$ satisfies $x_i = x_i^2$, whereas such equality is not guaranteed for points in $P$. This yields to the definition of the convex set
$$
\M{P} = \left\lbrace \begin{pmatrix} 1&x^\top \cr x&X \cr
\end{pmatrix} \succeq 0\ \Big|\ \eqref{lin:constr1}, \eqref{lin:constr2}, \eqref{lin:constr3}\ \text{hold} \right\rbrace,
$$
where \M{} is referred here as \textit{lift operator}.

\subsubsection*{Step 3 (Projection)} The projection of \M{P} onto the original $x$-space is
$$
\N{P} = \left\lbrace x \in  \mathbb{R}^n\ \Big|\ \exists\ Y\in M_+(P) \text{ with }\begin{pmatrix} 1 \\ x \end{pmatrix} = \text{diag}(Y) \right\rbrace.
$$
\citet{LoSc1991} showed that \N{P} is contained in $P$ and is a valid relaxation for $P_I$, i.e. $P_I \subseteq \N{P} \subseteq P$ holds. For the ease of presentation, in the reminder of this paper we will refer to both \M{P} and \N{P} as a relaxation of $P_I$, despite the former \say{lives} in the lifted variables space. Furthermore, note that all the linear inequalities $Ax \leq b$ are trivially satisfied by \add{points in} \M{P}.
\begin{remark}
    \label{rmk:inequality}
    Let $(a^k)^\top x \leq b^k$, for $k \in [m]$, be an inequality in the linear system defining $P$, then from the linear inequalities defining \M{P} one can sum~\eqref{lin:constr1} and~\eqref{lin:constr2}, for any $i \in [n]$, to recover the original constraint.
\end{remark}
\noindent
 Moreover, \citet{BoSaGu2017} noted that the definition of \M{P} and \N{P} does not depend on the representation of $P$ but only on the set of points in $P$. Hence, the following property holds\remove{:}
 
\begin{property}
	\label{pr:dominate}
	Given two polytopes $P, P^{\prime}$ with $P^{\prime} \subseteq P$ then $\M{P^{\prime}} \subseteq \M{P}$.
\end{property}

\noindent
Additionally, given any polytope $P$ containing $\ON{STAB}(G)$, \citet{LoSc1991} provided necessary conditions for an inequality $a^\top x \le b$ to be valid for \N{P}.

\begin{lemma}
	\label{th:contraction}
	 Let $a^\top x \le b$ be a valid inequality for $\ON{STAB}(G)$. W.l.o.g. assume $a \in \R{n}_+$ and $b \in \R{}_+$. If $a^\top x \le b$ is valid for $P \cap \{x:x_i=1\}$ for all $i \in V$, then $a^\top x \le b$ is valid for $N_+(P)$.
\end{lemma}

\section{Continuous relaxations of the SSP}
\label{sec:review}
In this section we recall known results about linear and semidefinite relaxation of the SSP. These are widely researched and we will only concentrate on elements relevant to our study. 

\subsection{Linear relaxations of the SSP}
\label{sec:linear}
Several classes of valid inequalities for $\ON{STAB}(G)$ have been identified since early seventies, \add{as is illustrated in}~\cite{Bo1998, GrLoSc2012}. Among them, {\em clique} inequalities $\sum_{i \in C} x_i \leq 1$, where $C \subseteq V$ is a maximal clique (set of pairwise adjacent vertices), induce facets of STAB$(G)$~\cite{Pa1973}. As customary in the literature, we define the polytope
$$
\qstab{G} = \{x \in [0,1]^{|V|}: \sum_{i \in C} x_i \leq 1, \; C \in {\cal K}\}
$$
where ${\cal K}$ denotes the collection of all maximal cliques of $G$. As clique inequalities clearly imply edge inequalities, STAB$(G) \subseteq$ \qstab{G} $\subseteq$ FRAC$(G)$ and inclusion is generally strict~\citep{GrLoSc2012}.  Clique inequalities \remove{in $G$} can be \replace{exponential in number}{exponentially many}, and the associated separation problem is NP-hard in the strong sense. However, effective separation heuristics make this class of inequalities computationally easy to handle, as shown in~\citet{MaRoSm2019}. 
\add{We say that a clique $C \in \mathcal{K}$ \textit{covers} an edge $\cbra{i, j} \in E$ if
$i, j \in C$.}
\replace{The}{An \emph{edge clique cover} of $G$ is a collection $\mathcal{C}$ of (not necessarily maximal) cliques 
of $G$ such that each edge in $E$ is covered by at least one clique in $\mathcal{C}$. Let ${\cal Q}$ be the set of all edge clique covers of $G$. Then, the} 0-1 linear program
\begin{align}
\max            & \sum\limits_{i \in V} \remove{w_i} x_i	  \nonumber \\
\mbox{s.t.}     & \sum\limits_{i \in C} x_i  \le \, 1   \quad \remove{(}C \in {\mathcal C}\remove{)} \label{eq:cov}\\
                & x_i                       \in \{0,1\} \quad \remove{(}i \in V\remove{)} \nonumber.
\end{align}
is a valid formulation for the SSP
\replace{if ${\mathcal C} \subseteq {\cal K}$ is any collection of maximal cliques covering all the edges of $G$ (i.e., both the endpoints of every edge are contained in at least one clique of $\mathcal{C}$)}{, for all ${\mathcal C} \in {\cal Q}$.} Extending the standard notation, we denote by \qstabc{G} the polytope associated with its continuous relaxation. \add{Unless stated otherwise, the results presented in the remainder of the paper are valid for any ${\mathcal C} \in {\cal Q}$.} 
It\replace{'s}{ is} not hard to see that $ \qstab{G} \subseteq \qstabc{G} \subseteq \ON{FRAC}(G)$. 
\replace{\qstabc{G} has a number of clique inequalities bounded by $|E|$}{
If ${\mathcal C} \in {\cal Q}$ is a \emph{minimal} edge clique cover, then $|\mathcal{C}| \leq |E|$. 
}

\remove{and typically contains considerably fewer inequalities than $\ON{FRAC(G)}$. In fact,}  \add{In practice,} \replace{F}{f}ormulations based on greedily computed 
\replace{minimum size collection}{ %
edge clique covers} $\mathcal{C}$, \add{ which are made up of maximal cliques, contain considerably fewer inequalities than those describing $\ON{FRAC(G)}$, and also provide a good approximation of $\qstab{G}$, as shown in several studies (see, e.g.,~\citet{LeRoSm2020} and references therein).}

Other {\em rank inequalities} of the form $\sum_{i \in S} x_i \le \alpha(G[S])$ have been studied for vertex subsets $S \subseteq V$ inducing special graphs (see e.g., \citet{Bo1998} and \citet{GiLeRoSm2009} for a survey). A well-known example arises when $S=H$, $|H| \ge 5$, induces a chordless cycle of odd cardinality, yielding the {\em odd-hole} inequality $\sum_{i \in H}x_i \le \lfloor |H|/2 \rfloor$.
Among other subgraphs, we recall {\em webs} and {\em antiwebs} introduced by~\citet{Tr1975}. Let $p$ and $q$ be integers satisfying $p > 2q + 1$ and $q > 1$,
and use arithmetic modulo $p$. A $(p,q)$−web is a graph with vertex set $\{1,...,p\}$ and with edges from $i$ to $\{i+q,\ldots,i−q\}$, for every $1 \leq i \leq p$. A $(p,q)$−antiweb is the complement of a $(p,q)$−web. The web inequalities take the form $\sum_{i \in W} x_i \le q$ for every vertex set $W$ inducing a $(p, q)$-web, and the antiweb inequalities take the form $\sum_{i \in AW} x_i \add{\leq} \lfloor p/q \rfloor$ for every vertex set $AW$ inducing a $(p,q)$-antiweb.

A different modeling rationale has been pursued by~\citet{DeTa1994} who considered the {\em surrogate} relaxation of formulation \add{based on}~\eqref{eq:edge} obtained by summing up, for every $i \in V$, the edge inequalities over all $j \in \Gamma(i)$. This yields the so called {\em nodal inequality} $\sum_{j \in \Gamma(i)} x_j + |\Gamma(i)| x_i \, \le \, |\Gamma(i)| $. \citet{MuCh1997} observed that this relaxation can be strengthened by replacing $|\Gamma(i)|$ with any value $r_i \ge \alpha(G[\Gamma(i)])$, which returns a class of 0-1 linear programs of the form

\begin{align}
\max		    &\,\, \sum\limits_{i \in V} \remove{w_i} x_i		 \nonumber \\
\mbox{s.t.}  	& \sum\limits_{j \in \Gamma(i)} x_j + \, r_i \, x_i \, \le \, r_i \quad \remove{(}i \in V\remove{)} \label{eq:nodal}\\
                & x_i \in \{0,1\}			\quad  \remove{(}i \in V\remove{)} \nonumber.
\end{align}

\noindent
The strongest formulation, corresponding to $r_i = \alpha(G[\Gamma(i)])$, has been extensively investigated in~\cite{LeRoSm2020}. Of course, computing coefficients $r_i = \alpha(G[\Gamma(i)])$ is NP-hard, as it amounts to solve the SSP on the subgraphs $G[\Gamma(i)]$. Nevertheless, \citet{LeRoSm2020} showed that these formulations can be handled efficiently for many graphs of interest. 
Here, we are interested in a further choice of coefficients $r_i$, namely, $r_i = \theta(G[\Gamma(i)])$ where $\theta(G[\Gamma(i)])$ is the \citet{Lo1979}~$\theta$ number of the  graph $G[\Gamma(i)]$. This represents an upper bound to $\alpha(G[\Gamma(i)])$ which can be computed in polynomial time up to an arbitrary precision, as reviewed in Section~\ref{sec:theta}. Then, we introduce the following notation
$$\nod{G}{r} = \cbra{x \in [0,1]^{|V|}: \sum_{j \in \Gamma(i)} x_j + \, r_i  \, x_i \, \le \, r_i, \hskip5pt i \in V }$$
\noindent
 to denote the polytope defined from~\eqref{eq:nodal} (and non-negativity), where $r \in\{\Gamma, \theta, \alpha\}$ specifies the coefficient used. Clearly, we have 
 $$\nod{G}{\alpha} \subseteq \nod{G}{\theta} \subseteq \nod{G}{\Gamma}.$$
It is worth noting that nodal inequalities in general do not imply edge inequalities, even in the case $r=\alpha$. On the other hand, if $r=\alpha$, they imply the \emph{wheel inequalities} (see~\citet{Bo1998} for their definition) which are not implied by clique inequalities.

\subsection{SDP relaxations of the SSP}

Three well-known SDP relaxations are related to our study.

\subsubsection{The Lov\'asz theta relaxation} 
\label{sec:theta} 

The Lov\'asz theta relaxation~\citep{Lo1979} has been introduced in 1979, yet we refer to the presentation of~\citet{LoSc1991}. Introduce the quadratic variable $X_{ij}$, representing the product $x_i x_j$ for all $\{i,j\} \in V \times V$ and let $X=xx^\top$ be the associated symmetric matrix of order $n = |V|$. An upper bound \replace{to}{on} $\alpha(G)$ is given by
\begin{align}
 \theta(G\remove{,w})=\max &\; \sum_{i \in V} \remove{w_i} x_i  \nonumber \\
         \mbox{s.t. }  &\; X_{ii} = x_i,  \quad	 i \in V, \label{eq:theta_sdp} \tag{SDP-$\theta$} \\
                      &\; X_{ij} = 0,    \quad  \{i,j\} \in E, \nonumber \\
             & \begin{pmatrix} 1&x^\top \cr x&X\cr  \end{pmatrix} \succeq 0.    \nonumber
\end{align}
\noindent
\remove{We simply denote it by $\theta(G)$ in the unweighted case.} 
The projection $\ON{TH}(G)$ of the feasible region of~\eqref{eq:theta_sdp} onto the subspace of the $x_i$ variables is known as the \emph{theta body}, a convex but not \add{a} polyhedral set in general. Remarkably, \citet{GrLoSc2012} proved that $\ON{STAB}(G) \subseteq \ON{TH}(G) \subseteq \qstab{G}$, where equality holds if and only if $G$ is perfect. This implies that there exists a polynomial-time separation algorithm for a class of inequalities which includes all clique inequalities. A thorough comparison between the original formulation in~\cite{Lo1979} and~\eqref{eq:theta_sdp} is provided by~\citet{GaLe2017}.
$\theta(G)$ can be computed in polynomial time up to an arbitrary precision and often provides a strong bound to the stability number $\alpha(G)$, typically better than those obtained from linear relaxations see, e.g.,~\cite{GiLeRoSm2009, GiRoSm2013, Ju1982, YiFa2006}. Another interesting feature of~\eqref{eq:theta_sdp} concerns its computational behavior. In fact,~\eqref{eq:theta_sdp} is handled better than similarly sized unstructured SDPs by general SDP solvers. Moreover, several algorithms have been designed to solve it efficiently as in, e.g.,~\cite{MaPoRe2009, PoReWi2006, GiLeRoSm2013}. 

Therefore, $\theta(G)$ achieves a good compromise between strength of the upper bound and computational burden. 
From a theoretical viewpoint, \citet{BuPa2006} showed that, unless P=NP,  no polynomially computable upper bound to $\alpha(G)$ which is provably smaller than $\theta(G)$ can be found. 
\replace{In this context}{However}, several stronger relaxations have been obtained by adding linear inequalities to \eqref{eq:theta_sdp}, as reviewed in the following subsections.

\subsubsection{The Schrijver relaxation}
\citet{Sc1979} observed that the PSD condition \add{in}~\eqref{eq:theta_sdp} does not imply the non-negativity on $X$. Hence, by adding the inequalities
\begin{align}
    \refstepcounter{equation} \tag{SDP-$\theta^+$a}\label{eq:theta_p_nonneg}
    X_{ij} \geq 0, \quad  \{i,j\} \in \bar{E}
\end{align}

\begingroup
\renewcommand\theequation{SDP-$\theta^+$}
\addtocounter{equation}{-1}
\refstepcounter{equation}\label{eq:theta_p_sdp}
\endgroup
\noindent
to~\eqref{eq:theta_sdp}, one yields a model that we denote as~\eqref{eq:theta_p_sdp} and the corresponding upper bound as $\theta^+(G)$.  Accordingly, $\ON{TH}^+(G)$ denotes the convex set obtained by its projection onto the subspace of $x_i$ variables.

\subsubsection{The Lov\'asz-Schrijver relaxation}

The Lov\'asz-Schrijver relaxation for the SSP is obtained by applying the operator \N{} to $\ON{FRAC}(G)$. 
\Nfrac{G} corresponds to~\eqref{eq:theta_sdp} plus the following linear inequalities:
\begin{subequations}
    \label{eq:LS}
    \begin{alignat}{3}
 X_{ij} &\geq 0, 						   &\quad& \{i,j\} \in \bar{E}, \\
 X_{ik} + X_{jk} &\leq x_k, 				   &\quad& \{i,j\} \in E, k \neq i,j, \label{eq:LS1}\\
 x_i + x_j + x_k &\leq 1 + X_{ik} + X_{jk}, &\quad& \{i,j\} \in E, k \neq i,j. \label{eq:LS2}
\end{alignat}
\end{subequations}
\noindent
\citet{LoSc1991} assessed the theoretical strength of \Nfrac{G} by showing that it satisfies several well-known classes of valid inequalities for $\ON{STAB}(G)$. Namely, they proved that \Nfrac{G} satisfies all clique, odd cycle, odd antihole and odd wheel inequalities.  \citet{GiLe2006} showed that, in addition, it satisfies all web inequalities. Optimizing $\sum_{i \in V} x_i$ over \Nfrac{G} returns an upper bound \replace{to}{on} $\alpha(G)$ that we denote as $\lambda(G)$.
From a more general perspective, it is interesting to observe the key role of the positive semidefiniteness constraint. In fact, \citet{BuVa2006} compared $\lambda(G)$ with the optimal value of the corresponding LP relaxation obtained by the Sherali-Adams procedure, observing a substantial improvement due to it.

Handling such a strong relaxation in practice is not straightforward. The first related computational results, obtained by a lift-and-project cutting plane method, appear in~\cite{BaCeCoPa1994}. Results on \replace{optimising}{optimizing} directly over \Mfrac{G} are illustrated in~\citet{Da2001} and, more recently, in~\citet{BuVa2006} using a \replace{specialised}{specialized} augmented Lagrangian method. The resulting upper bounds \replace{to}{on} $\alpha(G)$ may be considerably better than $\theta^+(G)$ but at the expense of a large increase in running times. In practice, this relaxation turns out to be hardly tractable for graphs of medium/large size. 

\subsubsection{Further SDP relaxations} 
\noindent
Variants of the previous relaxations are presented in~\citet{DuRe2007} and in~\citet{GrRe2003}. In the former, extra triangle inequalities are added to \Nfrac{G} to obtain a very strong relaxation. In the latter, relaxations are specialized with respect to graph density to improve efficiency. We refer to~\citet{GaLe2017} for an exhaustive review of these methods.
Two SDP-related methodologies are also presented in the literature. In~\citet{GaRe2020}\add{,} strong SDP relaxations of several combinatorial optimization problems are investigated by exploiting the \emph{exact subgraph inequalities}. Despite the power of this technique, the experience reveals that strengthening the natural SDP relaxation of the SSP is typically harder than other problems, e.g., \emph{max-cut}, and confirms that improving the $\theta$ bound towards $\alpha(G)$, even by a small amount, is actually challenging. Another effective method (presented in the context of the clique number) based on adding non-valid inequalities has been developed in~\citet{Lo2015}. This turns out to achieve quite strong bounds for both structured and random graphs with up to 150 vertices.

Summarizing, the evidence from state-of-the-art methodologies is that achieving even a small improvement of $\theta(G)$ requires a considerable additional computational cost. Namely, adding extra linear inequalities in the SDP models yields much harder SDPs and does require specialized methods. In the next section we introduce new relaxations with the goal of bridging such a gap.

\section{New SDP relaxations}
\label{sec:newSDP}

We now apply the lift operator \M{} to \qstabc{G} and \nod{G}{r} and compare the resulting SDP relaxations. Let us first observe that all of them contain variable bound constraints, the lifting step of which consists in multiplying $0 \le x_i \le 1$ by $x_j$ and $(1 − x_j)$, for $i, j \in V, i < j$:
\begin{subequations}
\begin{align}
x_i x_j \rightarrow  & \quad X_{ij} \ge 0,   \label{eq:liftbound1} \\
x_i (1 - x_j) \rightarrow & \quad x_i \ge X_{ij},  \label{eq:liftbound2} \\
(1 - x_i) x_j \rightarrow   & \quad x_j \ge X_{ij},  \label{eq:liftbound3} \\
(1-x_i)(1-x_j) \rightarrow & \quad x_i + x_j - X_{ij} \le 1. \label{eq:liftbound4}
\end{align}
\end{subequations}
\noindent
The linear inequalities \eqref{eq:liftbound1}--\eqref{eq:liftbound4} are known as~\citet{Mc1976} inequalities and describe the convex hull of the set $\{(x_i, x_j, y_{ij}) \in \{0,1\}^3 : y_{ij}=x_ix_j\}$. It is easy to see that they are not implied by the PSD condition~\eqref{eq:Schur} and should be explicitly added to the relaxations illustrated below unless implied by other constraints. 

\subsection{Relaxation \MqstabC{G}}
In order to investigate \MqstabC{G} we first describe the application of the operator \M{} to \qstab{G}. The lifting step consists in multiplying each clique inequality $\sum_{i \in C} x_i \leq 1$, $C \in {\cal K}$,  by $x_j$ and $(1-x_j)$, for all $j \in V$. Then, adding the condition $X_{ii} = x_i$, one obtains\remove{:}

\begin{subequations}
\begin{numcases}{x_j \cdot \left( \sum_{i \in C} x_i  \leq 1\right) \rightarrow}
\label{clq1} \sum_{i \in C \setminus \{j\}} X_{ij} \leq 0, &   if $j \in C$, \\
\label{clq2} \sum_{i \in C} X_{ij} - x_j \leq 0, &  if $j \not \in C$,
\end{numcases}
\begin{numcases}{(1 - x_j) \cdot \left( \sum_{i \in C} x_i  \leq 1\right) \rightarrow}
	\label{clq3}\sum\limits_{i \in C}x_i - \sum\limits_{i \in C \setminus \{j\}}X_{ij} \leq 1, &   if $j \in C$, \\
	\label{clq4}\sum\limits_{i \in C}(x_i - X_{ij}) + x_j \leq 1, &  if $j \not \in C$.
\end{numcases}
\end{subequations}
\noindent
Inequalities (\ref{clq2}), called \emph{clique-variable} inequalities, have been introduced in~\cite{GiLeRoSm2009}. One can observe that 
\begin{enumerate}
\item Constraints \eqref{clq1} imply $X_{ij} = 0$ for all $\{i,j\} \in E$ and constraints \eqref{clq3} reduce to the corresponding clique inequality $\sum_{i \in C} x_i \le 1, C \in {\cal K}$, which is already implied from~\eqref{clq2} and~\eqref{clq4} by Remark~\ref{rmk:inequality};

\item $X_{ij} = 0$ for all $\{i, j\} \in E$, $X_{ij} \ge 0$ for all $\{i,j\} \in \bar{E}$ and constraints \eqref{clq2} imply inequalities \eqref{eq:liftbound2} and \eqref{eq:liftbound3}; 
\item Constraints \eqref{clq1}, \eqref{clq2} and \eqref{clq4} imply \eqref{eq:liftbound4}.
\end{enumerate}
Then, the description of \Mqstab{G} reduces to

\begin{align*}
	\begin{split}
		\Mqstab{G} = \left\{	\begin{pmatrix}
		1 & x^\top \\
		x & X
	\end{pmatrix} \succeq 0 :\right. & \\
	\sum\limits_{i \in C}X_{ij} - x_j &\leq 0,\quad C \in \mathcal{K}, j \in V \setminus C, \\
	\sum\limits_{i \in C}(x_i - X_{ij}) + x_j &\leq 1,\quad C \in \mathcal{K}, j \in V \setminus C,\\
	x_{i} & = X_{ii}, \quad i \in V,\\
	X_{ij} &= 0,\quad    \{i,j\} \in E,\\
	X_{ij} &\geq 0, \quad \{i,j\} \in {\bar E}   \\
 x_i &	\ge 0, \quad i \in V \bigg\}.
\end{split}
\end{align*}

\noindent
Moreover, one has\remove{:}
\begin{theorem} 
\label{th:Q2FRAC}
$\Mqstab{G} \subseteq \Mfrac{G}$.
\end{theorem} 

\if\journal1
\begin{proof}
\else
\noindent
\textit{Proof }
\fi
    Since $\qstab{G} \subseteq \ON{FRAC}(G)$, the statement follows from Property~\ref{pr:dominate}. \qed

\if\journal1
\end{proof}
\else
\bigskip
\fi

\noindent
We denote \add{the resulting upper bound on $\alpha(G)$} by 
\begin{equation*}
    \mu(G) = \max\cbra{\sum_{i \in V} x_i: \begin{pmatrix} 1&x^\top \cr x&X\cr  \end{pmatrix} \in \Mqstab{G}}
\end{equation*}
The following example shows that the containment \replace{$\Mqstab{G} \subseteq \Mfrac{G}$}{of Theorem~\ref{th:Q2FRAC}} can be strict.

\begin{example}
    \label{ex:awfrac}
    Let $G$ be the $(10, 3)$-antiweb graph of Figure~\ref{fig:antiw}. We have $\mu(G) = 3$. Differently, optimizing the cardinality objective function $\sum_{i \in V} x_i$ over \Nfrac{G} one obtains the optimal solution $x_F = (0.3106, \ldots, 0.3106)$ with value $z^*_F = 3.106$. Therefore, $x_F \in \Nfrac{G}$ lies outside \Nqstab{G}. 
\end{example}

\noindent
The example indeed suggests a general result. In fact, we can show that, unlike \Nfrac{G}, \add{the} relaxation \Nqstab{G} implies the class of antiweb inequalities.

\begin{theorem}
\label{th:qstab_aw}
    \Nqstab{G} implies antiweb inequalities\add{.}
\end{theorem}
\if\journal1
\begin{proof}
\else
\noindent
\textit{Proof }
\fi
    Let $G$ be the $(p, q)$-antiweb graph with $V = \cbra{1, \ldots, p}$. In what follows, we use arithmetic modulo $p$. Note that, by definition, any consecutive $q$ nodes in $G$ form a clique. By denoting $\alpha = \left\lfloor p/q \right\rfloor$, the $(p, q)$-antiweb inequality is 
    $$ \sum\limits_{i \in V}x_i \leq \alpha.$$
    Given any $i \in V$, we denote \replace{with}{by} $G \setminus i$ the subgraph induced by the vertex set 
    $$V(G\setminus i) = V \setminus (\cbra{i} \cup \Gamma(i)) = \cbra{i + q,\ i + q + 1,\dots,\ i - q - 1,\ i - q}.$$ 
    The goal is to prove that inequality
    \begin{equation}
        \sum\limits_{j \in V(G \setminus i)}x_j \leq \alpha - 1, \label{eq:aw_proof}
    \end{equation}
    is valid for $\qstab{G \setminus i}$, for all $i \in V$, then the theorem's statement will follow from Lemma~\ref{th:contraction}. Let $i$ be any node in $V$ and consider the following $\alpha - 1$ cliques in $G \setminus i$
    \begin{align*}
        C^1 &= \cbra{i + q,\ i + q + 1,\dots,\ i + 2q - 1}, \\
        C^2 &= \cbra{i + 2q,\ i + 2q + 1,\dots,\ i + 3q - 1}, \\
            &\vdots \\
        C^k &= \cbra{i + kq,\ i + kq + 1,\dots, i + (k + 1)q - 1},\\
            &\vdots \\
        C^{\alpha - 1} &= \cbra{i + (\alpha - 1)q,\ i + (\alpha - 1)q + 1,\dots,\ i - q}. 
    \end{align*}
    Note that the set $\cbra{C^1, \dots, C^{\alpha - 1}}$ forms a partition of $V(G \setminus i)$, i.e., \remove{$\bigcup_{k = 1,\dots,\alpha-1} C^k = V(G \setminus i)$} $$\add{\bigcup_{k = 1,\dots,\alpha-1} C^k = V(G \setminus i)}$$ and $C^\ell \cap C^k = \emptyset$, for all $\ell, k =1,\dots, \alpha - 1$, with $\ell \not = k$. 
    Hence, the sum of clique inequalities $\sum_{j \in C^k}x_j \leq 1$, for $k=1,\dots,\alpha - 1$ \replace{yield}{yields} the inequality
    $$ \sum\limits_{i + q \leq j \leq i - q}x_j = \sum\limits_{j \in V(G \setminus i)}x_j \leq \alpha - 1, $$
    which is valid for $\ON{QSTAB}(G \setminus i)$. By symmetry, this argument is independent of the choice of $i$. \qed
\if\journal1
\end{proof}
\else
\bigskip
\fi

Theorem \ref{th:Q2FRAC} completes the theoretical picture for \Nqstab{G} along the line drawn by~\citet{LoSc1991}. In fact, all of the traditional combinatorial inequalities, namely, clique, odd cycle, odd antihole, odd wheel, web and antiweb inequalities turn out to be satisfied by \Nqstab{G}. It is interesting to mention that \citet{GiLeRoSm2009} proved the same result for a non-compact linear relaxation obtained by the application of the $N(K, K)$ lifting operator by Lov\'asz and Schrijver to \qstab{G}. The comparison between $N(K,K)$ and \N{} was indeed suggested by \citet{LoSc1991} as an interesting open issue.

\begin{figure}[h]
    \centering
    \includegraphics[scale=0.6]{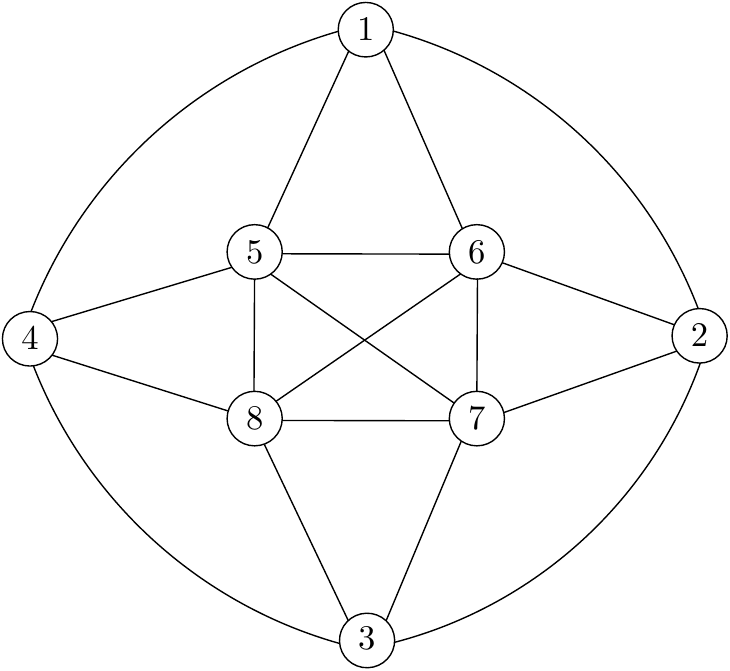}
    \caption{The $\hat G$ graph.}
    \label{fig:Ghat}
\end{figure}

Another insight of the progress of \Nqstab{G} towards \Nfrac{G} stems from the analysis of graph $\hat G$ of Figure~\ref{fig:Ghat} which has been shown  
in \cite{BiEsNaTu2014} that \Nfrac{\hat G} gives rise to a non-polyhedral relaxation. Here we have $\lambda(\hat G) = 2.146$. Now, optimizing over \Nqstab{\hat G} one gets \add{the stability number} $\mu(\hat G)= \alpha(\hat G) = 2$. \replace{This implies that \Nqstab{\hat G} is actually polyhedral.}
{This may suggest that \Nqstab{\hat G} is actually polyhedral.}

\noindent
From a computational perspective, a major issue of \Nqstab{G} is \add{that} the linear description of \qstab{G} \remove{that}, in general, contains exponentially many inequalities. \add{To get rid of this issue, we consider $\qstabc{G}$, for any edge clique cover $\cal C \in \cal Q$. The application of the lifting operator to this relaxation yields the convex set} \remove{However, as $\qstabc{G} \subseteq \ON{FRAC}(G)$, we still have $\MqstabC{G} \subseteq \Mfrac{G}$ and one can consider the convex set}
$$\MqstabC{G} = \cbra{ \begin{pmatrix}
    1 & x^\top \\
    x & X
\end{pmatrix} \succeq 0\, :\, x_i = X_{ii}, i\in V \mbox{ and } \eqref{clq1}-\eqref{clq4} \mbox{ hold} }.$$
\noindent
\add{Since $\qstabc{G} \subseteq \ON{FRAC}(G)$, we have $\MqstabC{G} \subseteq \Mfrac{G}$.}
\add{Also} note that \MqstabC{G} still implies McCormick inequalities \eqref{eq:liftbound1}--\eqref{eq:liftbound4}. 
Remarkably, it can be built in polynomial time as a \replace{suitable set}{minimal edge clique cover} ${\cal C} \add{ \in {\cal Q}}$ can be determined by a greedy algorithm (e.g., the one proposed in \citet{LeRoSm2020}), and the number of inequalities describing \qstabc{G} is \replace{bounded by}{at most} $|E|$ (for increasing graph density, it is also considerably smaller than $|E|$). Optimizing $\sum_{i \in V} x_i$ over \MqstabC{G} returns an upper bound that we denote as $\mu(G, {\cal C})$. We will experiment with such a relaxation in Section~\ref{sec:computation}.

\subsection{Relaxation \Nnod{G}{r}}
\label{subsec:m_nodal}

Let us now consider a generic nodal inequality, with $r_i \geq \alpha(G[\Gamma(i)])$, for all $i \in V$. The lifting step along with $x_i=X_{ii}$ returns:

$$	x_j \cdot \left( \sum_{\replace{h}{k} \in \Gamma(i)} x_{\replace{h}{k}} + r_i x_i \leq r_i \right)\rightarrow	$$
\begin{subequations}
	\begin{numcases} {\rightarrow}
		\label{nod1} \sum_{\replace{h}{k} \in \Gamma(i)} X_{i\replace{h}{k}} \leq 0, &   if $j =i$ \\
		\label{nod2} (1 - r_i)x_j + \sum_{\replace{h}{k} \in \Gamma(i) \setminus j} X_{j\replace{h}{k}} + r_i X_{ij} \leq 0, &  if $j \in \Gamma(i)$\\
		\label{nod3} \sum_{\replace{h}{k} \in \Gamma(i)} X_{j\replace{h}{k}} + r_i X_{ij} - r_i x_j \leq 0, &  if $j \not \in \Gamma(i)$
	\end{numcases}
\end{subequations}	

$$(1 - x_j) \cdot \left( \sum_{\replace{h}{k} \in \Gamma(i)} x_{\replace{h}{k}} + r_i x_i \leq r_i \right) \rightarrow$$
\begin{subequations}
	\begin{numcases}{\rightarrow}
		\label{nod4} \sum_{\replace{h}{k} \in \Gamma(i)} (x_{\replace{h}{k}} - X_{i\replace{h}{k}}) + r_i x_i\leq r_i, &  \hskip -10pt if $j =i$ \\
		\label{nod5} \sum_{\replace{h}{k} \in \Gamma(i) \setminus j} (x_{\replace{h}{k}} - X_{j\replace{h}{k}}) + r_i x_i + r_i x_j - r_i X_{ij} \leq r_i, & \hskip -10pt if $j \in \Gamma(i)$\\
		\label{nod6} \sum_{\replace{h}{k} \in \Gamma(i)} (x_{\replace{h}{k}} - X_{j\replace{h}{k}}) +r_i x_i + r_i x_j  -r_i X_{ij} \leq r_i, & \hskip -10pt if $j \not \in \Gamma(i)$
	\end{numcases}
\end{subequations}

\noindent
One can observe that:
\begin{enumerate}
	\item Constraints \eqref{nod1} imply $X_{ij} = 0$ for all $\{i,j\} \in E$; this is a relevant fact if one recalls that edge inequalities are not implied by nodal inequalities~\eqref{eq:nodal}; 
	\item $X_{ij} = 0$, for all $\{i, j\} \in E$, $X_{ij} \ge 0$ for all $\{i, j\} \in  \bar E$ and constraints \eqref{nod3} imply inequalities \eqref{eq:liftbound2} and \eqref{eq:liftbound3}; 
	\item $X_{ij} = 0$, for all $\{i, j\} \in E$, inequalities \eqref{eq:liftbound2} and \eqref{eq:liftbound3} and constraints \eqref{nod5} imply edge inequalities (and \eqref{eq:liftbound4} for $\{i,j\} \in E$);
    \item By $X_{ij} = 0$, for all $\{i, j\} \in E$,~\eqref{nod4} boils down to the original nodal inequality~\eqref{eq:nodal} and can be dropped, since Remark~\ref{rmk:inequality} can be applied to both pairs~\eqref{nod2},~\eqref{nod5} and~\eqref{nod3},~\eqref{nod6};
     \item $X_{ij} \ge 0$, for all $\{i, j\} \in \bar E$, inequalities \eqref{eq:liftbound2} and \eqref{eq:liftbound3} and constraints \eqref{nod6} imply \eqref{eq:liftbound4} for $\{i,j\} \in \bar E$.
 \end{enumerate}
\noindent
Thanks to these properties \Mnod{G}{r} boils down to

\begin{align*}
	\begin{split}
		\Mnod{G}{r} = \left\{	\begin{pmatrix}
			1 & x^\top \\
			x & X
		\end{pmatrix} \succeq 0 :\right. & \\
		 (1 - r_i) x_j  + \sum_{\replace{h}{k} \in \Gamma(i) \setminus j} X_{j\replace{h}{k}} \leq 0,&\quad i \in V, j \in \Gamma(i),   \\
		\sum_{\replace{h}{k} \in \Gamma(i)} X_{j\replace{h}{k}} + r_i X_{ij} - r_i x_j \leq 0,& \quad i \in V, j \not \in \Gamma(i), \\
		 \sum_{\replace{h}{k} \in \Gamma(i) \setminus j} (x_{\replace{h}{k}} - X_{j\replace{h}{k}}) + r_i x_i + r_i x_j \leq r_i, & \quad i \in V, j \in \Gamma(i),\\
		 \sum_{\replace{h}{k} \in \Gamma(i)} (x_{\replace{h}{k}} - X_{j\replace{h}{k}}) +r_i x_i + r_i x_j  -r_i X_{ij} \leq r_i, & \quad i \in V, j \not \in \Gamma(i),\\
		x_{i}  = X_{ii}, &\quad i \in V,\\
		X_{ij} = 0, & \quad    \{i,j\} \in E,                              \label{eq:edgezero} \\
		X_{ij} \geq 0, & \quad \{i,j\} \in {\bar E},   \\
		x_i 	\ge 0, & \quad i \in V \bigg\}.
	\end{split}
\end{align*}

\noindent
Then, \Mnod{G}{r} can be obtained by adding linear inequalities \eqref{nod2}, \eqref{nod3}, \eqref{nod5}, \eqref{nod6} and $X_{ij} \ge 0$, for $\{i,j\} \in \bar E$, to~\eqref{eq:theta_p_sdp} (it is straightforward to observe that inequalities \eqref{nod2} are implied by Mc-Cormick inequalities for $r = \Gamma$). Therefore, its projection \Nnod{G}{r} onto the $x$-space is at least as strong as the Schrijver relaxation\replace{:}{.}

\begin{theorem} 
\label{th:theta}
$\Nnod{G}{r} \subseteq \ON{TH}^+(G)$.
\end{theorem}
\noindent
This result also implies

\begin{corollary} 
    \Nnod{G}{r} implies all clique inequalities.
\end{corollary}

\begin{corollary} 
$\nu(G, r) = \max\cbra{\sum_{i \in V}x_i: \begin{pmatrix} 1&x^\top \cr x&X\cr  \end{pmatrix} \in \Mnod{G}{r}}\leq \theta^+(G)$.
\end{corollary} 

\noindent
Theorem \ref{th:theta} holds in particular for $r_i = |\Gamma(i)|$ corresponding to the case where the initial formulation is the surrogate relaxation of $\ON{FRAC(G)}$. Even in this case, applying the \N{} operator can recover the initial formulation's weakness and yield an upper bound that is at least as good as Schrijver's bound.  Furthermore, a hierarchy of relaxations derives from strengthening the lifting coefficient $r_i$, that is, looking at $r_i= \theta(G[\Gamma(i)])$ and $r_i= \alpha(G[\Gamma(i)])$. Thanks to Property~\ref{pr:dominate} we have

\begin{theorem}
    $ \Mnod{G}{\alpha} \subseteq \Mnod{G}{\theta} \subseteq \Mnod{G}{\Gamma}$.
\end{theorem}

\noindent
This hierarchy poses a complexity issue. The strongest relaxation comes from applying the operator to an LP relaxation, the construction of which, in general, cannot be carried out in polynomial time. An interesting polynomial special case arises when the subgraphs induced by neighbor sets are perfect and, according to the result in \cite{GrLoSc2012}, in this case, we have $\theta(G[\Gamma(i)]) = \alpha(G[\Gamma(i)])$, for all $i \in V$. Differently, $\Mnod{G}{\theta}$ can always be constructed (and solved) in polynomial time.
Section~\ref{sec:computation} documents the practical trade-off between strength and computational tractability for all relaxations in this hierarchy.

\subsection{Comparison among relaxations}
\label{sec:compare}
A natural question concerns with relationships between \Nfrac{G}, \allowbreak \Nqstab{G} \allowbreak and \allowbreak \Nnod{G}{r}. The following example shows that neither \Nnod{G}{r} contains \Nfrac{G} nor the reverse.

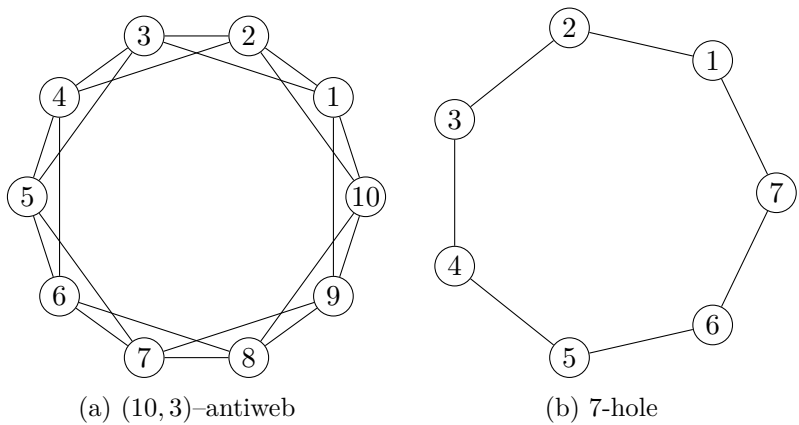
\begin{figure}[H]
    \centering
    \subcaptionbox{\label{fig:antiw}$(10,3)$--antiweb}{
		\tikzstyle{every node}=[circle, draw, fill=white!50,
		inner sep=0pt, minimum width=15pt]
		\begin{tikzpicture}[scale=.75]
			\tikzstyle{every node}=[circle, draw, fill=white!50,
			inner sep=0pt, minimum width=15pt]
			\def \n {10}
			\def \m {3}
			\def \radius {3cm}
			\def \margin {8} %
			\foreach \i in {1,...,10}
			{
				\ifthenelse{\i=4 \OR \i=3 \OR \i=1 \OR \i=2}{\node[draw] (n\i) at ({360/\n * (\i)}:\radius) {$\i$};}{\node[draw] (n\i) at ({360/\n * (\i)}:\radius) {$\i$};}
			}
			\foreach \i in {1,...,10}
			{
				\foreach \j in {1,...,2}
				{
					\pgfmathparse{int(Mod(\i + \j, \n)}
					\let \next \pgfmathresult
					\ifthenelse{\next=0}{\draw (n\i) -- (n\n);}{\draw (n\i) -- (n\next);}
				}	
			}
	\end{tikzpicture}}\quad
     \subcaptionbox{7-hole \label{fig:oddhole}}{
    		\tikzstyle{every node}=[circle, draw, fill=white!50,
    		inner sep=0pt, minimum width=15pt]
    			\begin{tikzpicture}[scale=.75]
    			\def \n {7}
    			\def \m {3}
    			\def \radius {3cm}
    			\def \margin {8} %
    			\foreach \i in {1,...,7}
    			{
    				\ifthenelse{\i=3 \OR \i=2 \OR \i=1}{\node[draw] (n\i) at ({360/\n * (\i)}:\radius) {$\i$};}{\node[draw] (n\i) at ({360/\n * (\i)}:\radius) {$\i$};}
    			}
    			\foreach \i in {1,...,7}
    			{
                    \pgfmathparse{int(Mod(\i + 1, \n)}
                    \let \next \pgfmathresult
                    \ifthenelse{\next=0}{\draw (n\i) -- (n\n);}{\draw (n\i) -- (n\next);}
    			}
	\end{tikzpicture}}
    \caption{Two example graphs.}
    \label{fig:example}
\end{figure}

\begin{example}
    Consider the odd-hole with seven vertices of Figure (\ref{fig:oddhole}). According to~\cite{LoSc1991}, \Nfrac{G} satisfies the facet-defining inequality $\replace{\sum_{i = 1}^8}{\sum_{i \in V}} x_i \le 3$ which implies $\lambda(G) = 3$. Differently, the optimal value from \Nnod{G}{r} is $\nu(G, \theta) = \nu(G, \alpha) = 3.317$ corresponding to solution $(0.47395, \ldots, 0.47395)$, which therefore falls outside \Nfrac{G}.     
    \noindent
    However, one can find a feasible solution to \Nfrac{G} which is outside \Nnod{G}{r} as well. Consider the $(10,3)$-antiweb graph $G$ of Figure~(\ref{fig:antiw}). Here, $\nu(G, \theta) = \nu(G, \alpha) = \alpha(G) = 3$. On the contrary, optimizing over \Nfrac{G} we obtain \add{the} solution $$x_{F} = (0.31055, \dots, 0,31055)$$ with value $3.105$ which is then outside \Nnod{G}{r}, for $r \in \{ \theta, \alpha \}$.
\end{example}

\begin{figure}
    \centering
    \includegraphics[scale=0.9]{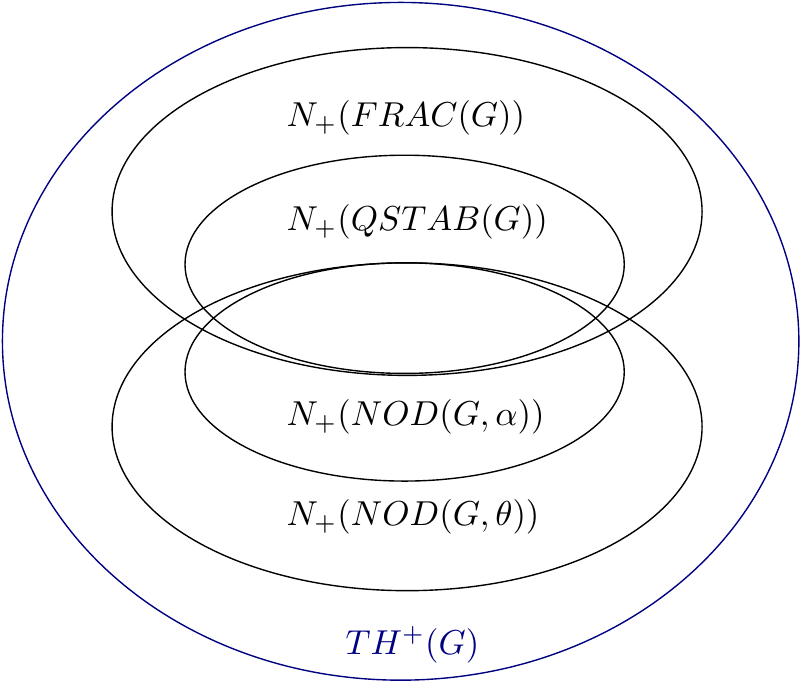}
    \caption{Containment relationships among relaxations.}
    \label{fig:relaxations}
\end{figure}

\noindent
This example also shows that the odd-hole inequalities are not implied by \Nnod{G}{r} even with $r=\alpha$. This appears as a weakness indicator of \Nnod{G}{r} if one recalls that odd-hole inequalities are contained in the first Chv\`atal closure of the 
\remove{edge} 
formulation \add{based on edge inequalities}~\eqref{eq:edge} \cite{Pa1973}. The overall theoretical picture is summarized in Figure \ref{fig:relaxations} 
\replace{, and the above results}{. The results related to the implication of odd-hole and antiweb inequalities} seem to suggest that \Nqstab{G} (as well as \Nfrac{G}) should be stronger than \Nnod{G}{r}. Interestingly, the computational analysis will reveal significant supremacy of the latter as graphs get dense, for $r \in \{ \theta, \alpha \}$.

\begin{table}
    \centering
    \begin{tabular}{l@{\hskip 0.35in}l@{\hskip 0.35in}}
\toprule
Optimal value & Relaxation \\
\midrule 
$\lambda(G)$ & \Nfrac{G} \\[.5ex]
$\mu(G)$ & \Nqstab{G} \\[.5ex]
$\mu(G, {\cal C})$ & \NqstabC{G} \\[.5ex]
$\nu(G, r)$ & \Nnod{G}{r}, for $r \in \{\Gamma, \theta, \alpha\}$ \\
\bottomrule
\end{tabular}

    \caption{Notation for SDPs optimal values.}
    \label{tab:sdp_relaxations}
\end{table}

\section{Implementation}
\label{sec:algorithms}
The software implemented for this work, along with the datasets used in the experiments (see Section~\ref{sec:computation}), is available at~\url{https://github.com/febattista/SDP_lift_and_project}. The repository comprises the following Python (\texttt{.py}) and MATLAB (\texttt{.m}) files:
\begin{itemize}
    \item \texttt{LinearFormulations.py}: this module creates all linear formulations presented in Section~\ref{sec:linear} and exports them in standard \texttt{.lp} format. It relies on the max-clique solver {\tt cliquer}~\cite{Os2002} and \texttt{ADMMsolver.py} (see below) to compute coefficients for nodal inequalities~\eqref{eq:nodal}.
    \item \texttt{SDPLifting.py}: starting from a linear formulation in standard \texttt{.lp} format, this module implements the \M{} operator described in Section~\ref{sec:lift_and_project}, exporting the resulting SDP in the \texttt{.mat} MATLAB format.
    \item \texttt{ADMMsolver.py}: this module implements the SDP solver proposed in~\cite{Ba2023, BaDe2021} and is exclusively employed as a subroutine during the LP formulation.
    \item \texttt{kelley\_cutting\_plane.m}: this module implements the cutting plane method described in Algorithm~\ref{alg:separation}.
\end{itemize}
The following Python scripts are provided to simplify the use of the software:
\begin{itemize}
    \item \texttt{install.py} downloads and builds all required dependencies (i.e. {\tt SDPNAL+}~\cite{SuToYuZh2020} and {\tt cliquer}~\cite{Os2002}).
    \item \texttt{parameters.py} sets up the computational experiments, allowing for (i) selection of the datasets to be used and (ii) tuning of all the parameters of Algorithm~\ref{alg:separation} and the SDP solvers involved.
    \item \texttt{analyze\_results.py}  collects the data from the experiments and replicates all the tables presented in this work.
\end{itemize}
\subsubsection*{Selection of the SDP solver} Under mild assumptions, semidefinite programs can be solved in polynomial time up to any arbitrary fixed precision~\citep{NeNe1994}. Intuitively, the computational burden mainly depends on the order $n$ of the matrix variable and the number $m$ of linear constraints. However, the density of the latter should also be taken into account as most SDP solvers rely on sparse matrices implementations.
\textit{Interior-point methods} (IPMs)~\cite{MOSEK, sedumi, sdpt3} provide a good-precision solution in reasonable time for \say{small} and \say{medium} SDPs (both in term of $n$ and $m$). However, they require to store and factorize large Hessian matrices at each iteration, which may become prohibitive for large SDPs.
In fact, we preliminarily tested the IPM implemented in the commercial state-of-the-art solver {\tt MOSEK}~\cite{MOSEK} and observed that problems with $m \approx 10^5$ and $n \approx 200 $ are intractable. Although this experiment already provided interesting insights on the strength of the relaxations, testing larger graphs is crucial for a reliable experimental picture. Therefore, we resorted to Alternating Direction Methods of Multipliers (ADMMs), which, as a variant of Augmented Lagrangian methods, is a popular first-order alternative to IPMs and can scale on much larger SDPs at the price of possibly slowing down the convergence to high-precision solutions. For a detailed overview, we refer the reader to~\cite{ Ba2023, BaDe2021, CeDeGaWi2021, PoReWi2006,  SuToYuZh2020, WeGoYi2010, WiZh2022}. Furthermore, ADMMs allow the efficient handling of bound constraints (e.g., nonnegativity) on the matrix variable via iterative projections on the cone of positive semidefinite matrices. In our case, it enables the computation of $\theta^+(G)$ with no substantial additional effort with respect to $\theta(G)$. We selected {\tt SDPNAL+} proposed in~\cite{SuToYuZh2020}, which represents a state-of-the-art ADMM general-purpose SDP solver implemented in MATLAB equipped with efficient subroutines in \texttt{C} included via \texttt{.mex} files, and has proven to be numerically stable in our experiments.

\subsubsection*{Model building} The \M{} operator has been applied to four different compact LP relaxations of the SSP: \qstabc{G} and \nod{G}{r} for $r \in \{\Gamma, \theta, \alpha \}$. The construction of these LP models entails different complexity levels.  \qstabc{G} requires \add{the computation of} an \replace{collection}{edge clique cover} \add{$\mathcal{C} \in \mathcal{Q}$} \remove{covering all the edges of $G$}.
\add{This task has been accomplished by the greedy algorithm presented in~\citet{LeRoSm2020} and is reported in Algorithm~\ref{alg:greedy_clq_cover}.
Starting from the empty collection $\mathcal{C} = \emptyset$, the heuristic iteratively adds a maximal clique of $G$ to $\mathcal{C}$, until every edge is covered. At each iteration, the selected maximal clique $C_k$ covers at least one new edge not yet covered by cliques in $\mathcal{C}$ (line~\ref{alg:clq:line1}).
Consequently, even if it can be shown that the resulting edge clique cover $\mathcal{C}$ may not be minimal, its size is at most $|E|$ and, in practice, comprises significantly fewer cliques~\cite{LeRoSm2020}.}
\remove{This is determined by the greedy heuristic described in \cite{LeRoSm2020}.}
\begin{algorithm}[ht]
    \caption{\add{Greedy edge clique cover}}
    \label{alg:greedy_clq_cover}
    \begin{algorithmic}[1]
        \State \add{\textbf{Input}: Graph $G = (V, E)$}
        \State \add{\textbf{Output}: An edge clique cover $\mathcal{C}$ of $G$}
        \State \add{Initialize $E^\prime \gets E,$ $G^\prime \gets (V, E^\prime)$ $\mathcal{C} \gets \emptyset,\ k \gets 0$}
        \While{\add{$E^\prime \not= \emptyset$}} 
        \State \add{$k \gets k + 1$}
        \State \add{$u \gets \arg\max_{i \in V}\cbra{|\Gamma_{G^\prime}(i)|}$}
        \State \add{$v \gets \arg\max_{i \in \Gamma_{G^\prime}(u)}\cbra{|\Gamma_{G^\prime}(i)|}$ }
        \Comment{\add{$E^\prime \not= \emptyset \Rightarrow\Gamma_{G^\prime}(u)$ is non-empty}}
        \State \add{Find a maximal clique $C_k$ of $G$ such that $u, v \in C_k$} \label{alg:clq:line1}
        \State \add{$\mathcal{C} \gets \mathcal{C} \cup \cbra{C_k}$}
        \State \add{$E^\prime \gets E^\prime \setminus E(C_k)$}
        \EndWhile
    \end{algorithmic}
\end{algorithm}

\noindent
Building \nod{G}{\Gamma} is straightforward,  as $r_i = |\Gamma(i)|$. On the contrary, \nod{G}{\theta} requires the evaluation of $r_i = \theta(G[\Gamma(i)])$, for all $i \in V$.  This corresponds to solve $|V|$ SDPs, namely, those corresponding to \eqref{eq:theta_sdp} on the subgraph induced by $\Gamma(i)$. In practice, the SDPs are completely decomposed and can be solved in parallel.  Moreover, $r_i$ can be tightened to $\left \lfloor \theta(G[\Gamma(i)]) \right\rfloor$ and one can stop the computation as long as the SDP solver returns a valid upper (dual) bound on $\theta(G[\Gamma(i)])$. Several approaches to compute a valid dual bound throughout solver's iterations and a posteriori have been proposed for ADMM algorithms (see, e.g.,~\cite{CeDeGaWi2021, JaChKe2008, WiZh2022}). Their specific application in SDP relaxations of the SSP have been implemented in~\cite{Ba2023, BaDe2021}. Then, for each iteration of the ADMM implemented in \texttt{ADMMsolver.py}, a valid dual bound is computed and this, in turn, allows stopping the computation of $\theta$ as soon as a \say{moderate} precision is achieved. In Table~\ref{tab:dimacs_coeff}, we report the average/max/min CPU times to calculate coefficients $\theta(G[\Gamma(i)])$ \replace{on a subset of the instances used in our computational experiments}{on a subset of the instances used in our computational experiments, whose detailed introduction is deferred to Section~\ref{sec:instances}}.
 Finally, building \nod{G}{\alpha} is the most challenging case as it requires calculating $\alpha(G[\Gamma(i)])$, for all $ \replace{I}{i} \in V$. Despite the theoretical hardness of determining $\alpha(G[\Gamma(i)])$,~\citet{LeRoSm2020} showed that in practice, the computational burden to build \nod{G}{\alpha} is, in fact, accessible for large classes of graphs. We refer the interested reader to their discussion; here we mention that even in this case the evaluation of coefficients $r_i = \alpha(G[\Gamma(i)])$ can be parallelized and, as in the previous case, we report in Table~\ref{tab:dimacs_coeff} the average/max/min CPU times to evaluate $\alpha(G[\Gamma(i)])$ via the max-clique solver {\tt cliquer}~\citep{Os2002}.

\begin{sidewaystable}
    \centering
\if\journal1
    \vspace{300pt}
\fi
    \scalebox{.90}{\begin{tabular}{ll|rrrrrrrrrrrrrrrr}
& Graph & 
\rotatebox{80}{\texttt{brock800\_1}} & 
\rotatebox{80}{\texttt{brock800\_2}} & 
\rotatebox{80}{\texttt{brock800\_3}} & 
\rotatebox{80}{\texttt{brock800\_4}} & 
\rotatebox{80}{\texttt{C500-9}} & 
\rotatebox{80}{\texttt{p\_hat300-1}} & 
\rotatebox{80}{\texttt{p\_hat300-2}} & 
\rotatebox{80}{\texttt{p\_hat300-3}} & 
\rotatebox{80}{\texttt{p\_hat500-1}} & 
\rotatebox{80}{\texttt{p\_hat500-2}} & 
\rotatebox{80}{\texttt{p\_hat500-3}} & 
\rotatebox{80}{\texttt{p\_hat700-1}} & 
\rotatebox{80}{\texttt{p\_hat700-2}} & 
\rotatebox{80}{\texttt{p\_hat700-3}} & 
\rotatebox{80}{\texttt{sanr400\_0.5}} & 
\rotatebox{80}{\texttt{sanr400\_0.7}} \\
\toprule
\multirow[c]{3}{*}{$\theta_i$}
& Min & 2.94 & 2.83 & 3.07 & 2.81 & 0.76 & 1.96 & 1.39 & 0.44 & 4.63 & 2.61 & 0.58 & 8.89 & 5.30 & 0.77 & 1.27 & 0.69 \\
& Max & 4.98 & 4.97 & 4.94 & 4.91 & 10.92 & 7.43 & 20.35 & 9.33 & 10.73 & 78.22 & 13.28 & 25.81 & 123.65 & 24.58 & 2.64 & 2.04 \\
& Average & 3.93 & 3.90 & 3.93 & 3.93 & 2.40 & 3.53 & 6.14 & 1.86 & 7.66 & 15.95 & 3.84 & 15.06 & 31.25 & 6.98 & 1.78 & 1.09 \\
\toprule
\multirow[c]{3}{*}{$\alpha_i$}
& Min & 0.89 & 1.08 & 0.68 & 0.63 & $<0.01$ & $<0.01$ & $<0.01$ & $<0.01$ & $<0.01$ & $<0.01$ & $<0.01$ & $<0.01$ & $<0.01$ & $<0.01$ & 0.01 & $<0.01$ \\
& Max & 28.87 & 34.41 & 26.17 & 25.91 & $<0.01$ & $<0.01$ & 0.02 & 0.02 & 0.02 & 1.83 & 2.74 & 0.07 & 62.41 & 108.48 & 0.04 & 0.13 \\
& Average & 6.00 & 6.17 & 5.93 & 6.04 & $<0.01$ & $<0.01$ & $<0.01$ & $<0.01$ & $<0.01$ & 0.10 & 0.10 & 0.02 & 1.64 & 2.78 & 0.02 & 0.03 \\
\bottomrule
\end{tabular}
}
    \caption{CPU time statistics for computing coefficients $\theta$ and $\alpha$ for \eqref{eq:nodal} on DIMACS graphs.}
    \label{tab:dimacs_coeff}
\end{sidewaystable}

\subsubsection*{Lifting operation} The function \texttt{m\_plus\_lifting()} in \texttt{SDPLifting.py} carries out Steps 1 and 2 of the procedure outlined in Section~\ref{sec:lift_and_project}. Initially, it parses the polytope $P$, provided as input through a standard \texttt{.lp} format text file. Then, it constructs the symmetric matrices corresponding to constraints~\eqref{lin:constr1}, \eqref{lin:constr2} and \eqref{lin:constr3}. Finally, the model is exported in \texttt{.mat} file format, ready to be loaded into MATLAB and solved using \texttt{SDPNAL+}.

\subsubsection*{SDP optimization} In general, we provide as input to \texttt{SDPNAL+} SDPs of the form 

$$ \ON{opt}(\ON{R}) = \max_{Y \in \ON{R}} \langle C, Y \rangle,$$
where
$$ \ON{R} = \cbra{Y \in {\cal S}^+_{n+1}\ :\ Y\geq 0, \langle A_i, Y \rangle \leq b_i,\ i \in [m]},$$
and $m \in \mathbb{N}$, $C, A_1, \ldots, A_m \in {\cal S}_{n+1}$, $b \in \R{m}$. 
In a preliminary experiment we observed that \texttt{SDPNAL+} converges quickly on~\eqref{eq:theta_p_sdp} while it requires very large (and practically unacceptable) convergence times on  \MqstabC{G} and \Mnod{G}{r} relaxations. This is due to the number of linear constraints in the former, and their density in the latter. However, we proved that \MqstabC{G} and \Mnod{G}{r} are contained in~\eqref{eq:theta_p_sdp}, thus one can resort to a classical cutting plane method to optimize over them.
The cutting plane algorithm starts optimizing over~\eqref{eq:theta_p_sdp}.
At the $k$-th iteration, given the optimal solution $Y^k \in {\cal S}^+_{n+1}$, one can compute
$$ V_i(Y^k) = \langle A_i, Y^k \rangle - b_i,\ i \in [m].$$
If $ V_i(Y^k) > \varepsilon$, for a \say{small} $\varepsilon > 0$, then the $i$-th constraint is violated by the current solution. Let $c \le m$ be a positive integer parameter, and let $I_c(Y^k) = \cbra{i_1, i_2, \ldots, i_c}$ denote the set of indexes corresponding to the  $c$ most violated constraints. Then, the constraints $\cbra{Y \in {\cal S}_{n+1}\ :\ \langle A_i, Y \rangle \leq b_i,\ i \in I_c(Y^k)}$ are added to the current SDP relaxation and a new optimal solution $Y^{k+1}$ is computed. The algorithm terminates either when no violated constraints are found or when a tailing off condition is met. 
The detailed algorithm is reported in Algorithm~\ref{alg:separation}.

\begin{algorithm}[ht]
    \caption{Cutting-plane scheme}
    \label{alg:separation}
    \begin{algorithmic}[1]
        \State \textbf{Input}: Any SDP $\ON{R} \subseteq \eqref{eq:theta_p_sdp}$, $c \in \mathbb{N}$, $\varepsilon > 0, \delta \in \R{}_+$
        \State Initialize $\Pi^0 \gets {\cal S}_{n+1}, \Delta^0 = +\infty,\ k=0$ 
        \State Compute $p^0 = \ON{opt}(\mbox{\ref{eq:theta_p_sdp}})$ and let $Y^0$ be its optimal solution \label{alg:line1}
        \While{$Y^k$ violates constraints in $\ON{R}$ \textbf{or} $\Delta^k > \delta$} \label{alg:line5}
        \State Compute the vector $V(Y^k)$ of violations and $I_c(Y^k)$ \label{alg:line2}
        \State $ \Pi^{k+1}\gets \Pi^k \cap \cbra{Y \in {\cal S}_{n+1}\ :\ \langle A_i, Y \rangle \leq b_i,\ i \in I_c(Y^k)}$
        \State Compute $p^{k+1} = \ON{opt}(\mbox{\ref{eq:theta_p_sdp}} \cap\ \Pi^{k+1})$ and let $Y^{k+1}$ be its optimal solution \label{alg:line3}
        \State $\Delta^{k+1} \gets |p^{k + 1} - p^k|$ \label{alg:line4}
        \State $k \gets k + 1$
        \EndWhile
    \end{algorithmic}
\end{algorithm}

\section{Computational Experiments}
\label{sec:computation}

We now illustrate the results obtained by solving relaxations \MqstabC{G} and \Mnod{G}{r} for $r \in \{\Gamma, \theta, \alpha \}$
with Algorithm \ref{alg:separation}. Experiments are designed to address three main questions:
\begin{itemize}
    \item[(i)] Can relaxations based on operator $M_+(\cdot)$ provide upper bounds on $\alpha(G)$ that are significantly stronger than $\theta(G)$?
    \item[(ii)] What is the additional computational cost for achieving such an improvement?
    \item[(iii)]  How much does the choice of $r$ affect the quality of the bound from \Mnod{G}{r}?
\end{itemize}
As a byproduct, we document the quality of the bound provided by \allowbreak \Mnod{G}{\theta}, \allowbreak that is, the strongest bound which can be computed in polynomial time in our hierarchy. 
\add{Finally, we investigate the extent to which the bound $\mu(G, \mathcal{C})$ is sensitive to the choice of $\mathcal{C}$ and, in particular, the degree to which it can be improved with a more refined selection of the edge clique cover.}

Experiments are run on a computer with Intel(R) Xeon(R) CPU E5-2698 v4 @ 2.20GHz, 256GB RAM, under OS Ubuntu 16.04.7 LTS (GNU/Linux 4.15.0-128-generic x86\_64). After preliminary tuning experiments, the parameters of Algorithm~\ref{alg:separation} have been set to $\replace{\epsilon}{\varepsilon} = 10^{-3}, \delta = 10^{-1}$ and $c = 1000$ with 7200 seconds time limit (excluding formulation building time). The precision of the SDP solver \texttt{SDPNAL+} is set to $10^{-6}$. The precision of \texttt{ADMMsolver.py} for evaluating coefficients $\theta_i$ is set to $10^{-4}$.
 Denoting by $\ON{opt}(\ON{R})$ the optimal (with the given precision) value of relaxation $\ON{R}$ computed by Algorithm \ref{alg:separation} at termination, tables below report, for each relaxation, four statistics: the percentage gap to $\alpha(G)$, that is, $\frac{\ON{opt}(\ON{R})- \alpha(G)}{\alpha(G)} \%$; the number of cutting plane iterations and the corresponding CPU time; the number of added cuts. The latter is detailed for each class of inequalities and a class is not reported if none of its members has ever been detected as violated.

\subsection{Instances}
\label{sec:instances}
The numerical experiments are based on the following two collections of graphs, available at~\url{https://github.com/febattista/SDP_lift_and_project} in the standard edge-list format.

\subsubsection*{Random graphs} A Erd\"os–R\'enyi graph $G(n, p)$ has $n$ vertices and each edge appears with probability $p \in [0, 1]$. We have considered the collection of graphs from \cite{LeRoSm2020}, generated by 
$$n \in \{200, 250, 300\},\quad p \in \{0.1, 0.2, 0.3, 0.4, 0.5, 0.6, 0.7, 0.8, 0.9\},$$
where, for each combination of $n$ and $p$, five instances with different random seeds have been created for a total of 135 graphs. 

\subsubsection*{DIMACS graphs} Graphs from the Second DIMACS Implementation Challenge~\citep{JoTr1996} form the standard benchmark for max-clique algorithms. We complemented the graphs to convert the max-clique instances into SSP instances  (\texttt{DSJC*} graphs are not complemented, as these belong to the ‘‘coloring” benchmark set). We excluded the \texttt{c-fat}, \texttt{johnson}, and \texttt{san} graphs because the integrality gap is completely closed by linear relaxations, so they do not provide information to our analysis. We also dropped \texttt{keller5} and the largest \texttt{p-hat} instances (with 1000 and 1500 vertices) as the SDP relaxations turned out to be too large to be handled. The final collection includes 38 graphs. Statistics about the time required to compute lifting coefficients of linear formulations \Nnod{G}{\theta} and \Nnod{G}{\alpha} are reported in Table \ref{tab:dimacs_coeff}. Recall that their computation can be carried out in parallel, so the impact on the overall time is often negligible.

\subsection{Experiment 1: random graphs}
Table \ref{tab:random1} collects gap statistics for Erd\"os–R\'enyi graphs, where the first column reports the integrality gap of the LP relaxation \nod{G}{\alpha}.
\replace{Tables \ref{tab:random_times} and \ref{tab:random_cuts} report cutting plane iterations, CPU times, and, respectively, the number of added cuts for each class of inequalities.}{Tables~\ref{tab:random_times} and \ref{tab:random_cuts} report cutting plane iterations and CPU times, and the number of added cuts for each class of inequalities, respectively.} 

Let us first notice that $\lambda(G) \simeq \theta^+(G)$ in all cases but $n=200$, $p \in \{0.1, 0.9\}$, evidencing a rather disappointing behavior of \Mfrac{G}. On the contrary, replacing edge inequalities with a compact collection of clique inequalities in the initial formulation pays off, as $\mu(G, {\cal C})$ is always smaller than $\lambda(G)$. The average portion of additional gap closed by $\mu(G, {\cal C})$ w.r.t. $\lambda(G)$ is $6.7\%$. This is achieved by a few cutting plane iterations, adding about one to three thousand inequalities of class \eqref{clq2}, resulting in an average $80\%$ CPU time increase. Interestingly, no inequalities in the class \eqref{clq4} have been detected as violated.

This experiment clearly classifies relaxations strength as a function of edge probability $p$ and graph size $n$. When $p \le 0.4$,  \MqstabC{G} turns out to be the strongest relaxation and the sparser the graph, the larger the gap closed by $\mu(G, {C})$ w.r.t. other relaxations, while formulations based on nodal inequalities fail to improve on $\theta^+(G)$. However, the advantage of \MqstabC{G} over other relaxations tends to thin out as $n$ increases.  The outcome is completely reversed when $p \ge 0.5$, where relaxations based on nodal inequalities are more competitive and, mainly, \Mnod{G}{\alpha} is the clear winner. $\nu(G, \alpha)$ closes $45.3\%$ of the gap left by $\lambda(G)$ and $40.0\%$ of the one left by $\mu(G, {\cal C})$, on average. To our knowledge, such an improvement on $\theta^+(G)$ has rarely been documented in the literature on SDP approaches, even for smaller graphs. Another nice evidence is that this advantage of \Mnod{G}{\alpha} is independent of the graph size. Observe that, also in this case, one class of inequalities is predominant, namely, \eqref{nod2}. Notice also that some violated inequalities in \eqref{nod5} and \eqref{nod6} have only been detected for the densest graphs $p=0.9$. Very good upper bounds are also achieved by \Mnod{G}{\theta} when $p \ge 0.7$. In fact, $\nu(G, \theta)$ remarkably improves on $\theta^+(G)$ and, unlike $\nu(G, \alpha)$, can be computed in polynomial time. Also, the computational effort required to achieve strong bounds is reasonable. Differently, the weakest relaxation \Mnod{G}{\Gamma} always returns $\nu(G, \Gamma)= \theta^+(G)$, that is, never closes additional gap in our tests. Finally, notice that for highly dense graphs, the gap of the LP relaxation \nod{G}{\alpha} may be smaller than that of the weakest SDP relaxations. However, it is largely improved by the strongest SDP relaxation, which confirms the operator's power and the PSD constraint's key role.

\begin{table}
    \centering
    	\setlength{\tabcolsep}{0.2em}

\begin{tabular}{ll@{\hskip 0.20in}r@{\hskip 0.20in}r@{\hskip 0.20in}r@{\hskip 0.20in}r@{\hskip 0.20in}r@{\hskip 0.20in}r@{\hskip 0.20in}r}
\toprule
 \multicolumn{2}{l}{Graph} &  \multicolumn{7}{c}{\% gap} \\
 \toprule
$n$ & $p$ & $\ON{NOD}(G, \alpha)$ & $\theta^+(G)$ & $\lambda(G)$ & $\mu(G, \mathcal{C})$ & $\nu(G, \Gamma)$ & $\nu(G, \theta)$ & $\nu(G, \alpha)$ \\

\midrule
\multirow[c]{9}{*}{200}
 & 0.1 &  46.957 & 19.330 & 18.723 & \textbf{16.033} & 19.330 & 19.330 & 19.330 \\
 & 0.2 & 102.539 & 27.924 & 27.912 & \textbf{26.060} & 27.924 & 27.924 & 27.924 \\
 & 0.3 &  93.334 & 30.774 & 30.773 & \textbf{29.393} & 30.774 & 30.774 & 30.774 \\
 & 0.4 &  73.543 & 30.550 & 30.550 & \textbf{29.322} & 30.550 & 30.550 & 30.537 \\
 & 0.5 &  56.527 & 30.597 & 30.597 & \textbf{29.355} & 30.597 & 30.593 & 29.540 \\
 & 0.6 &  39.889 & 25.708 & 25.708 & 24.318 & 25.708 & 25.143 & \textbf{20.138} \\
 & 0.7 &  30.543 & 25.579 & 25.579 & 23.679 & 25.579 & 19.897 & \textbf{12.494} \\
 & 0.8 &   9.720 & 13.474 & 13.474 & 10.819 & 13.474 & 2.871 & \textbf{1.836} \\
 & 0.9 &   9.600 & 10.123 & 10.120 & 5.346 & 10.123 & \textbf{5.057} & \textbf{5.057} \\
 \midrule
\multirow[c]{9}{*}{250} 
 & 0.1 &  64.521 & 27.371 & 27.112 & \textbf{24.908} & 27.371 & 27.371 & 27.371 \\
 & 0.2 & 118.348 & 35.462 & 35.459 & \textbf{34.284} & 35.462 & 35.462 & 35.462 \\
 & 0.3 & 103.026 & 39.136 & 39.136 & \textbf{38.370} & 39.136 & 39.136 & 39.136 \\
 & 0.4 &  78.520 & 38.335 & 38.335 & \textbf{37.737} & 38.335 & 38.335 & 38.312 \\
 & 0.5 &  54.333 & 33.563 & 33.563 & 32.943 & 33.563 & 33.563 & \textbf{31.421} \\
 & 0.6 &  46.178 & 39.416 & 39.416 & 38.669 & 39.416 & 39.163 & \textbf{28.532} \\
 & 0.7 &  39.114 & 38.531 & 38.531 & 37.459 & 38.531 & 36.305 & \textbf{23.679} \\
 & 0.8 &  19.867 & 20.012 & 20.012 & 18.473 & 20.012 & 9.419 & \textbf{4.433} \\
 & 0.9 &   5.950 & 11.717 & 11.716 & 9.524 & 11.717 & \textbf{3.219} & \textbf{3.219} \\
 \midrule
\multirow[c]{9}{*}{300} 
 & 0.1 &  75.661 & 31.124 & 31.047 & \textbf{29.389} & 31.124 & 31.124 & 31.124 \\
 & 0.2 & 137.919 & 46.238 & 46.237 & \textbf{45.551} & 46.238 & 46.238 & 46.238 \\
 & 0.3 & 113.088 & 48.426 & 48.426 & \textbf{48.033} & 48.426 & 48.426 & 48.426 \\
 & 0.4 &  89.040 & 50.228 & 50.228 & \textbf{49.918} & 50.228 & 50.228 & 50.166 \\
 & 0.5 &  62.433 & 46.061 & 46.061 & 45.791 & 46.061 & 46.061 & \textbf{41.439} \\
 & 0.6 &  42.224 & 39.956 & 39.956 & 39.595 & 39.956 & 39.848 & \textbf{26.444} \\
 & 0.7 &  23.650 & 31.619 & 31.619 & 31.059 & 31.619 & 28.272 & \textbf{11.347} \\
 & 0.8 &  23.333 & 29.809 & 29.809 & 28.718 & 29.809 & 19.807 & \textbf{11.667} \\
 & 0.9 &   7.090 & 23.160 & 23.160 & 20.318 & 23.160 & 5.512 & \textbf{5.486} \\
 \bottomrule
\end{tabular}

    \caption{$N_+(\cdot)$ bounds on random instances. The \% gap of the LP relaxation $\nod{G}{\alpha}$ is reported as a reference. The best gaps are shown in boldface.}
    \label{tab:random1}
\end{table}

\begin{table}
    \centering
	\setlength{\tabcolsep}{0.2em}

    \begin{tabular}{ll@{\hskip 1em}r@{\hskip 1em}rr@{\hskip 1em}rr@{\hskip 1em}rr@{\hskip 1em}rr@{\hskip 1em}rr}
\toprule
 \multicolumn{2}{c@{\hskip 1em}}{Graph} &  \multicolumn{1}{c@{\hskip 1em}}{$\theta^+(G)$} & \multicolumn{2}{c@{\hskip 1em}}{$\lambda(G)$} & \multicolumn{2}{c@{\hskip 1em}}{$\mu(G, \mathcal{C})$} & \multicolumn{2}{c@{\hskip 1em}}{$\nu(G, \Gamma)$} & \multicolumn{2}{c@{\hskip 1em}}{$\nu(G, \theta)$} & \multicolumn{2}{c}{$\nu(G, \alpha)$} \\
 \toprule
$n$ & $p$ & Time & \# iter & Time & \# iter & Time & \# iter & Time & \# iter & Time & \# iter & Time \\
\midrule
\multirow[c]{9}{*}{200} 
 & 0.1 & 7.24 & \replace{1.6}{2.0} & 26.64 & 3.2 & 48.24 & 0.2 & 8.74 & 0.2 & 8.64   & 0.2 & 8.76 \\
 & 0.2 & 4.29 & 1.0 & 9.10  & 2.0 & 13.96 & 0.0 & 4.34 & 0.0 & 4.34   & 0.0 & 4.34 \\
 & 0.3 & 2.64 & 0.8 & 4.43  & 2.0 & 8.53  & 0.0 & 2.71 & 0.0 & 2.71   & 0.0 & 2.71 \\
 & 0.4 & 1.74 & 0.2 & 2.65  & 2.0 & 6.51  & 0.0 & 1.82 & 0.2 & 2.13   & 1.0 & 3.77 \\
 & 0.5 & 1.64 & 0.0 & 2.46  & 2.0 & 6.37  & 0.0 & 1.74 & 1.0 & 3.60   & \replace{1.0}{1.8} & 12.12 \\
 & 0.6 & 1.82 & 0.0 & 2.85  & 2.0 & 7.00  & 0.0 & 1.93 & 1.0 & 6.56   & 2.0 & 22.02 \\
 & 0.7 & 2.18 & 0.0 & 3.26  & 2.0 & 10.90 & 0.0 & 2.30 & 2.0 & 29.37  & 2.4 & 58.03 \\
 & 0.8 & 3.05 & 0.0 & 4.42  & 2.0 & 14.14 & 0.0 & 3.18 & 1.4 & 49.62  & 1.4 & 40.79 \\
 & 0.9 & 5.43 & 0.6 & 10.33 & 2.0 & 37.73 & 0.0 & 5.57 & 2.6 & 236.08 & 2.6 & 236.57 \\
\midrule
\multirow[c]{9}{*}{250}
 & 0.1 & 9.62 & \replace{1.0}{1.6} & 28.92 & 3.0 & 56.76 & 0.0 & 9.68 & 0.0 & 9.68   & 0.0 & 9.68 \\
 & 0.2 & 5.15 & 1.0 & 11.53 & 2.0 & 18.05 & 0.0 & 5.25 & 0.0 & 5.25   & 0.0 & 5.25 \\
 & 0.3 & 4.05 & 0.4 & 5.82  & 2.0 & 11.34 & 0.0 & 4.19 & 0.0 & 4.18   & 0.0 & 4.18 \\
 & 0.4 & 2.40 & 0.2 & 4.17  & 1.0 & 5.12  & 0.0 & 2.57 & 0.0 & 2.57   & 1.0 & 4.95 \\
 & 0.5 & 2.10 & 0.0 & 4.08  & 1.0 & 4.78  & 0.0 & 2.30 & 0.6 & 3.51   & 2.0 & 23.23 \\
 & 0.6 & 2.05 & 0.0 & 4.21  & 1.0 & 4.79  & 0.0 & 2.26 & 1.0 & 6.46   & 3.0 & 53.63 \\
 & 0.7 & 2.07 & 0.0 & 4.66  & 1.0 & 5.87  & 0.0 & 2.30 & 2.0 & 29.51  & 3.0 & 65.98 \\
 & 0.8 & 2.52 & 0.0 & 5.04  & 1.0 & 8.39  & 0.0 & 2.80 & 2.8 & 112.29 & 2.6 & 162.95 \\
 & 0.9 & 5.77 & 0.2 & 9.94  & 1.2 & 20.54 & 0.0 & 6.05 & 3.0 & 430.18 & 3.0 & 429.29 \\
\midrule 
\multirow[c]{9}{*}{300} 
 & 0.1 & 10.92 & 1.0 & 25.16 & 3.0 & 68.62 & 0.0 & 11.04 & 0.0 & 11.03  & 0.0 & 11.03 \\
 & 0.2 & 8.02  & 0.6 & 13.87 & 2.0 & 25.53 & 0.0 & 8.18  & 0.0 & 8.19   & 0.0 & 8.18 \\
 & 0.3 & 5.26  & 0.0 & 6.88  & 1.0 & 9.61  & 0.0 & 5.47  & 0.0 & 5.47   & 0.0 & 5.48 \\
 & 0.4 & 3.31  & 0.0 & 5.72  & 1.0 & 6.59  & 0.0 & 3.60  & 0.0 & 3.58   & 1.0 & 7.26 \\
 & 0.5 & 2.67  & 0.0 & 5.62  & 1.0 & 5.64  & 0.0 & 3.01  & 0.0 & 3.00   & 2.0 & 31.62 \\
 & 0.6 & 2.54  & 0.0 & 6.24  & 1.0 & 5.89  & 0.0 & 2.90  & 1.0 & 7.51   & 3.0 & 66.31 \\
 & 0.7 & 2.37  & 0.0 & 6.86  & 1.0 & 5.80  & 0.0 & 2.84  & 2.0 & 40.36  & 4.0 & 135.51 \\
 & 0.8 & 3.78  & 0.0 & 7.42  & 1.0 & 10.76 & 0.0 & 4.31  & 2.2 & 97.62  & 3.6 & 197.18 \\
 & 0.9 & 6.61  & 0.0 & 11.98 & 1.8 & 30.75 & 0.0 & 7.11  & 4.8 & 1250.1 & 5.0 & 1566.21 \\
\bottomrule
\end{tabular}

    \caption{CPU time and cutting plane's iterations on random graphs.}
    \label{tab:random_times}
\end{table}

\begin{table}
    \centering
    	\setlength{\tabcolsep}{0.2em}
    \begin{tabular}{ll@{\hskip 2.0em}r@{\hskip 2.0em}r@{\hskip 2.0em}rrrr@{\hskip 2.0em}rrrr}
	
	\toprule
	\multicolumn{2}{c@{\hskip 2.0em}}{Graph} & \multicolumn{10}{c}{Number of violated inequalities added in Algorithm~\ref{alg:separation}} \\
	\toprule
	\multicolumn{2}{c@{\hskip 2.0em}}{} & \multicolumn{1}{c@{\hskip 2.0em}}{$\lambda(G)$} & \multicolumn{1}{c@{\hskip 2.0em}}{$\mu(G, \mathcal{C})$} & \multicolumn{4}{c@{\hskip 2.0em}}{$\nu(G, \theta)$}  & \multicolumn{4}{c}{$\nu(G, \alpha)$} \\
$n$ & $p$ &  \eqref{eq:LS1} & \eqref{clq2} & \eqref{nod2}& \eqref{nod3} & \eqref{nod5}& \eqref{nod6}& \eqref{nod2}& \eqref{nod3} & \eqref{nod5}& \eqref{nod6}\\
 \midrule 
\multirow[c]{9}{*}{200} 
 & 0.1  & 1018 & 3126 & 0    & 0  & 0 & 0 & 0 & 0       & 0 & 0             \\
 & 0.2  & 30   & 2000 & 0    & 0  & 0 & 0 & 0 & 0       & 0 & 0             \\
 & 0.3  & 2    & 1517 & 0    & 0  & 0 & 0 & 0 & 0       & 0 & 0             \\
 & 0.4  & 0    & 1326 & 0    & 1  & 0 & 0 & 8 & 2       & 0 & 0             \\
 & 0.5  & 0    & 1232 & 7    & 1  & 0 & 0 & 569 & 65    & 0 & 0          \\
 & 0.6  & 0    & 1213 & 322  & 28 & 0 & 0 & 1459 & 128  & 0 & 0     \\
 & 0.7  & 0    & 1247 & 1508 & 52 & 0 & 0 & 2137 & 157  & 0 & 0    \\
 & 0.8  & 0    & 1273 & 1332 & 27 & 0 & 0 & 1395 & 5    & 0 & 0      \\
 & 0.9  & 1    & 1280 & 1800 & 0  & 567 & 60 & 1800 & 0 & 567 & 60 \\
 \midrule
\multirow[c]{9}{*}{250} 
 & 0.1 & 615 & 3000 & 0 & 0 & 0 & 0 & 0 & 0           & 0 & 0             \\
 & 0.2 & 7   & 1821 & 0 & 0 & 0 & 0 & 0 & 0           & 0 & 0             \\
 & 0.3 & 0   & 1281 & 0 & 0 & 0 & 0 & 0 & 0           & 0 & 0             \\
 & 0.4 & 0   & 972  & 0 & 0 & 0 & 0 & 21 & 3          & 0 & 0            \\
 & 0.5 & 0   & 907  & 1 & 1 & 0 & 0 & 1095 & 145      & 0 & 0        \\
 & 0.6 & 0   & 862  & 201 & 11  & 0 & 0 & 2799 & 201  & 0 & 0     \\
 & 0.7 & 0   & 908  & 1122 & 32 & 0 & 0 & 2924 & 76   & 0 & 0     \\
 & 0.8 & 0   & 994  & 2759 & 41 & 0 & 0 & 2378 & 122  & 0 & 0    \\
 & 0.9 & 0   & 755  & 2400 & 0  & 232 & 24 & 2400 & 0 & 232 & 24 \\
 \midrule
\multirow[c]{9}{*}{300} 
 & 0.1 & 256 & 3000 & 0 & 0 & 0 & 0 & 0 & 0           & 0 & 0              \\
 & 0.2 & 1   & 1462 & 0 & 0 & 0 & 0 & 0 & 0           & 0 & 0              \\
 & 0.3 & 0   & 946  & 0 & 0 & 0 & 0 & 0 & 0           & 0 & 0              \\
 & 0.4 & 0   & 657  & 0 & 0 & 0 & 0 & 69 & 8          & 0 & 0             \\
 & 0.5 & 0   & 534  & 0 & 0 & 0 & 0 & 1822 & 178      & 0 & 0         \\
 & 0.6 & 0   & 547  & 128 & 4   & 0 & 0 & 2828 & 172  & 0 & 0       \\
 & 0.7 & 0   & 662  & 1772 & 60 & 0 & 0 & 3939 & 61   & 0 & 0      \\
 & 0.8 & 0   & 848  & 2138 & 62 & 0 & 0 & 3574 & 26   & 0 & 0      \\
 & 0.9 & 0   & 1179 & 4000 & 0  & 460 & 63 & 4200 & 0 & 568 & 73  \\
 \bottomrule
\end{tabular}

    \caption{Number of added cutting planes on random instances.}
    \label{tab:random_cuts}
\end{table}

\subsection{Experiment 2: DIMACS instances}
The DIMACS test set, a standard for max-clique/stable set studies, includes, for the most part, structured graphs from various applications. The breadth of the graph collection supports the relevance of our analysis. To our knowledge, DIMACS graphs with more than 300-400 vertices have never been documented in SDP-based studies~\cite{DuRe2007, GaRe2020, GaSiWi2022, Lo2015}. The graphs we have considered are significantly larger than those in case studies applying the \N{} operator in~\cite{ BuVa2006, Da2001}. In our experience, experimenting with larger graphs is crucial to draw reliable conclusions about the strength of relaxations for the stable set problem, as it typically tends to degrade with graph size. Looking at Table \ref{tab:dimacs}, as already noted for random graphs, $\lambda(G)$ does not prove to be better than $\theta^+(G)$: some improvement is observed in two cases (\texttt{DSJC125-1} and \texttt{MANN\_a27}), it is negligible in nine cases and null in the remaining ones. Concerning the relaxations we have proposed, \Mnod{G}{\alpha} shows the most promising results. $\nu(G, \alpha)$ is the best bound in 11 cases while $\mu(G, {\cal C})$ is the winner in as many as 23 cases. Nevertheless, $\nu(G, \alpha)$ achieves the smallest gap on average, thanks to its remarkable strength on selected instances. According to the evidence gathered from Erd\"os–R\'enyi graphs, \MqstabC{G} is more suited for sparse instances, and \Mnod{G}{\alpha} turns out to be particularly tight as the density increases.  $\mu(G, {\cal C})$ is the best bound for the \texttt{C} and \texttt{MANN} collections and the sparsest among \texttt{p\_hat} graphs. On the other hand, $\nu(G, \alpha)$ is very strong for \texttt{DSJC500-5}, \texttt{p\_hat300-1}, \texttt{p\_hat500-1}, \texttt{p\_hat700-1}, \texttt{sanr400\_0.5}, where it closes more than $50\%$ of the gap left by
$\theta^+(G)$ and $\lambda(G)$ (about $74\%$ for \texttt{p\_hat-1} graphs). The computational price for this is reported in Table \ref{tab:dimacs_times}: computing $\mu(G, \alpha)$ is about 27 times slower than computing $\theta^+(G)$ except \texttt{p\_hat700-1} where it is 185 times slower.
But this experiment also shows that the tightness of \Mnod{G}{\alpha} is less sensitive to the graph size with respect to the other relaxations. This results in an interesting outcome for the hard \texttt{brock800} collection that, with a density of $35 \%$, one expects to be favorable to \Mfrac{G} and \MqstabC{G}. $\nu(G, \alpha)$ is the only bound stronger than $\theta^+(G)$ on these graphs, and the improvement is significant. Table \ref{tab:dimacs_times} shows that computing $\nu(G, \alpha)$ is about four times slower than computing $\theta^+(G)$ (including the time to build the LP relaxation), an acceptable price for such a valuable result. Table \ref{tab:dimacs_cuts} also reveals the key contribution to closing the gap given by inequalities \eqref{nod2}.

Also in this experiment, the polynomial bound $\nu(G, \theta)$ improves on $\theta^+(G)$ in six cases by a meaningful amount (namely, \texttt{DSJC125-5}, \texttt{MANN\_a27}, \texttt{p\_hat300-1}, \texttt{p\_hat500-1}, \texttt{p\_hat700-1}). The comparison between $\nu(G, \theta)$ and $\nu(G, \alpha)$ reveals that strengthening the lifting coefficient pays off, as in eleven cases, $\nu(G, \alpha)$ is significantly smaller than $\nu(G, \theta)$, in other eight cases the improvement is negligible while no improvement is observed in the remaining sixteen cases. On the other hand, $\nu(G, \Gamma)$ has never enhanced $\theta^+(G)$ even if a few violated cuts have occasionally been detected.

Looking at Tables \ref{tab:dimacs_times} and \ref{tab:dimacs_cuts}, one can observe that relevant improvements on $\theta^+(G)$ are achieved by adding several inequalities in a few iterations of the cutting plane algorithm. Similarly to what is observed for random graphs, inequalities \eqref{clq4} are only added in a few special graphs, namely the extremely sparse \texttt{MANN} graphs. Differently, classes of inequalities \eqref{nod2} and \eqref{nod3} are relevant to describe \Mnod{G}{\theta}, \Mnod{G}{\alpha}, while no violated members of classes \eqref{nod5} and \eqref{nod6} have ever been found.
\add{Interestingly, our numerical results indicate that inequalities generated during the lifting step by multiplying by $x_i \ge 0$ are more frequently violated than those obtained by multiplying by $1 - x_i \ge 0$. This trend is consistent across all reported relaxations, with only a few exceptions.}

Algorithm \ref{alg:separation} has shown a nice numerical behavior even for the largest graphs in the DIMACS collection, and re-optimization is cost-effective (Table \ref{tab:dimacs_times}): its extra time is only paid when the cuts yield a significant improvement of the upper bound. In contrast, in the other cases, CPU times do not deviate significantly from those of $\theta^+(G)$.

Overall, the computational analysis shows that the new SDP relaxations \MqstabC{G} and, mainly, \Mnod{G}{\alpha}, may be a viable option to compute upper bounds \replace{to}{on} $\alpha(G)$ which are remarkably stronger than $\theta(G)$ at a reasonable computational price.  Interestingly, some experimental evidence complements the theoretical insights illustrated in Section~\ref{sec:compare}, where \MqstabC{G} looks to be likely to achieve \replace{the tightest}{tighter} bounds than \Mnod{G}{\alpha}.

\begin{sidewaystable}
\if\journal1
    \vspace{390pt}
\fi
    \centering
    \scriptsize
    \begin{tabular}{lc@{\hskip 0.2in}rr@{\hskip 0.2in}rr@{\hskip 0.2in}rr@{\hskip 0.2in}rr@{\hskip 0.2in}rr@{\hskip 0.2in}rr}
\toprule
 \multicolumn{2}{c@{\hskip 0.2in}}{Graph} & \multicolumn{2}{c@{\hskip 0.2in}}{$\theta^+(G)$} & \multicolumn{2}{c@{\hskip 0.2in}}{$\lambda(G)$} & \multicolumn{2}{c@{\hskip 0.2in}}{$\mu(G, \mathcal{C})$} & \multicolumn{2}{c@{\hskip 0.2in}}{$\nu(G, \Gamma)$} & \multicolumn{2}{c@{\hskip 0.2in}}{$\nu(G, \theta)$} & \multicolumn{2}{c}{$\nu(G, \alpha)$}  \\
 Name & $d$(\%) & UB & \% gap & UB & \% gap & UB & \% gap & UB & \% gap & UB & \% gap & UB & \% gap \\
\midrule
\texttt{brock200\_1}  & 25 & 27.20 & 29.508 & 27.20 & 29.505 & \textbf{27.12} & \textbf{29.119} & 27.20 & 29.508 & 27.20 & 29.508 & 27.20 & 29.508 \\ 
\texttt{brock200\_2}  & 50 & 14.13 & 17.758 & 14.13 & 17.758 & 14.12 & 17.674 & 14.13 & 17.758 & 14.13 & 17.757 & \textbf{14.02} & \textbf{16.795} \\ 
\texttt{brock200\_3}  & 39 & 18.67 & 24.479 & 18.67 & 24.479 & \textbf{18.65} & \textbf{24.356} & 18.67 & 24.479 & 18.67 & 24.478 & 18.67 & 24.467 \\ 
\texttt{brock200\_4}  & 34 & 21.12 & 24.242 & 21.12 & 24.242 & \textbf{21.10} & \textbf{24.101} & 21.12 & 24.242 & 21.12 & 24.242 & 21.12 & 24.242 \\ 
\midrule
\texttt{brock400\_1}  & 25 & 39.33 & 45.670 & 39.33 & 45.670 & \textbf{39.33} & \textbf{45.662} & 39.33 & 45.670 & 39.33 & 45.670 & 39.33 & 45.670 \\ 
\texttt{brock400\_2}  & 25 & 39.20 & 35.160 & 39.20 & 35.160 & \textbf{39.20} & \textbf{35.156} & 39.20 & 35.160 & 39.20 & 35.160 & 39.20 & 35.160 \\ 
\texttt{brock400\_3}  & 25 & 39.16 & 26.324 & 39.16 & 26.324 & \textbf{39.16} & \textbf{26.310} & 39.16 & 26.324 & 39.16 & 26.324 & 39.16 & 26.324 \\ 
\texttt{brock400\_4}  & 25 & 39.23 & 18.883 & 39.23 & 18.883 & \textbf{39.23} & \textbf{18.877} & 39.23 & 18.883 & 39.23 & 18.883 & 39.23 & 18.883 \\ 
\midrule
\texttt{brock800\_1}  & 35 & 41.87 & 82.032 & 41.87 & 82.032 & 41.87 & 82.032 & 41.87 & 82.032 & 41.87 & 82.032 & \textbf{41.55} & \textbf{80.646} \\ 
\texttt{brock800\_2}  & 35 & 42.10 & 75.435 & 42.10 & 75.435 & 42.10 & 75.435 & 42.10 & 75.435 & 42.10 & 75.435 & \textbf{41.73} & \textbf{73.896} \\ 
\texttt{brock800\_3}  & 35 & 41.88 & 67.530 & 41.88 & 67.530 & 41.88 & 67.530 & 41.88 & 67.530 & 41.88 & 67.530 & \textbf{41.51} & \textbf{66.043} \\ 
\texttt{brock800\_4}  & 35 & 42.00 & 61.541 & 42.00 & 61.541 & 42.00 & 61.541 & 42.00 & 61.541 & 42.00 & 61.541 & \textbf{41.61} & \textbf{60.039} \\ 
\midrule
\texttt{C125-9}       & 10 & 37.55 & 10.431 & 36.75 & 8.074 & \textbf{36.23} & \textbf{6.566} & 37.55 & 10.431 & 37.54 & 10.424 & 37.55 & 10.429 \\ 
\texttt{C250-9}       & 10 & 55.82 & 26.856 & 55.71 & 26.618 & \textbf{55.20} & \textbf{25.446} & 55.82 & 26.856 & 55.82 & 26.856 & 55.82 & 26.856 \\ 
\texttt{C500-9}       & 10 & 83.58 & 46.630 & 83.58 & 46.630 & \textbf{83.51} & \textbf{46.515} & 83.58 & 46.630 & 83.58 & 46.630 & 83.58 & 46.630 \\ 
\midrule
\texttt{DSJC125.1}    &  9 & 38.04 & 11.896 & 37.28 & 9.650 & \textbf{36.79} & \textbf{8.207} & 38.04 & 11.895 & 38.04 & 11.881 & 38.04 & 11.881 \\ 
\texttt{DSJC125.5}    & 50 & 11.40 & 14.021 & 11.40 & 14.017 & 11.36 & 13.598 & 11.40 & 14.021 & 11.38 & 13.810 & \textbf{11.35} & \textbf{13.531} \\ 
\texttt{DSJC125.9}    & 90 & \textbf{4.00} & \textbf{0.000} & \textbf{4.00} & \textbf{0.000} & \textbf{4.00} & \textbf{0.000} & \textbf{4.00} & \textbf{0.000} & \textbf{4.00} & \textbf{0.000} & \textbf{4.00} & \textbf{0.000} \\ 
\texttt{DSJC500-5}    & 50 & 22.57 & 73.621 & 22.57 & 73.621 & 22.57 & 73.621 & 22.57 & 73.621 & 22.57 & 73.621 & \textbf{20.54} & \textbf{58.014} \\ 
\midrule
\texttt{MANN\_a9}     &  7 & 17.48 & 9.219 & 17.09 & 6.811 & \textbf{17.00} & \textbf{6.250} & 17.48 & 9.219 & 17.47 & 9.201 & 17.47 & 9.201 \\ 
\texttt{MANN\_a27}    &  1 & 132.76 & 5.367 & 131.11 & 4.057 & \textbf{131.01} & \textbf{3.980} & 132.76 & 5.367 & 131.99 & 4.751 & 131.99 & 4.751 \\ 
\texttt{johnson32-2-4}& 12 & \textbf{16.00} & \textbf{0.000} & \textbf{16.00} & \textbf{0.000} & \textbf{16.00} & \textbf{0.000} & \textbf{16.00} & \textbf{0.000} & \textbf{16.00} & \textbf{0.000} & \textbf{16.00} & \textbf{0.000} \\ 
\texttt{keller4}      & 35 & 13.47 & 22.417 & 13.46 & 22.388 & \textbf{13.39} & \textbf{21.709} & 13.47 & 22.417 & 13.45 & 22.236 & 13.45 & 22.236 \\ 
\midrule
\texttt{p\_hat300-1}  & 76 & 10.02 & 25.253 & 10.02 & 25.253 & 10.02 & 25.194 & 10.02 & 25.253 & 9.58 & 19.691 & \textbf{8.58} & \textbf{7.288} \\ 
\texttt{p\_hat300-2}  & 51 & 26.71 & 6.855 & 26.58 & 6.317 & \textbf{26.49} & \textbf{5.950} & 26.71 & 6.855 & 26.70 & 6.813 & 26.69 & 6.768 \\ 
\texttt{p\_hat300-3}  & 26 & 40.70 & 13.057 & 40.62 & 12.829 & \textbf{40.44} & \textbf{12.320} & 40.70 & 13.057 & 40.70 & 13.054 & 40.70 & 13.053 \\ 
\midrule
\texttt{p\_hat500-1}  & 75 & 13.01 & 44.533 & 13.01 & 44.533 & 13.01 & 44.526 & 13.01 & 44.533 & 12.65 & 40.583 & \textbf{10.68} & \textbf{18.721} \\ 
\texttt{p\_hat500-2}  & 50 & 38.56 & 48.306 & 38.44 & 47.863 & \textbf{38.38} & \textbf{47.598} & 38.56 & 48.306 & 38.56 & 48.293 & 38.54 & 48.225 \\ 
\texttt{p\_hat500-3}  & 25 & 57.81 & 15.622 & 57.76 & 15.524 & \textbf{57.67} & \textbf{15.342} & 57.81 & 15.622 & 57.81 & 15.622 & 57.81 & 15.621 \\ 
\midrule
\texttt{p\_hat700-1}  & 75 & 15.05 & 36.774 & 15.05 & 36.774 & 15.05 & 36.773 & 15.05 & 36.774 & 14.86 & 35.047 & \textbf{11.30} & \textbf{2.709} \\ 
\texttt{p\_hat700-2}  & 50 & 48.44 & 10.091 & 48.30 & 9.776 & \textbf{48.24} & \textbf{9.646} & 48.44 & 10.091 & 48.44 & 10.087 & 48.41 & 10.023 \\ 
\texttt{p\_hat700-3}  & 25 & 71.76 & 15.734 & 71.71 & 15.660 & \textbf{71.64} & \textbf{15.545} & 71.76 & 15.734 & 71.75 & 15.734 & 71.75 & 15.733 \\ 
\midrule
\texttt{sanr200\_0.7} & 30 & 23.63 & 31.296 & 23.63 & 31.291 & \textbf{23.59} & \textbf{31.048} & 23.63 & 31.296 & 23.63 & 31.296 & 23.63 & 31.296 \\ 
\texttt{sanr200\_0.9} & 10 & 48.90 & 16.440 & 48.63 & 15.786 & \textbf{48.03} & \textbf{14.356} & 48.90 & 16.440 & 48.90 & 16.440 & 48.90 & 16.440 \\ 
\texttt{sanr400\_0.5} & 50 & 20.18 & 55.217 & 20.18 & 55.217 & 20.18 & 55.216 & 20.18 & 55.217 & 20.18 & 55.217 & \textbf{19.10} & \textbf{46.951} \\ 
\texttt{sanr400\_0.7} & 30 & 33.97 & 61.746 & 33.97 & 61.746 & \textbf{33.97} & \textbf{61.741} & 33.97 & 61.746 & 33.97 & 61.746 & 33.97 & 61.746 \\ 
\bottomrule
\textbf{Means} &  &     - & 30.832 &     - & 30.528 &              - &          30.248 &     - & 30.832 &     - & 30.489 &     - & \textbf{27.772} \\
\bottomrule
\end{tabular}
    \caption{Lift-and-project bounds on DIMACS instances. The column $d(\%)$ \replace{report}{reports} the density of each graph.}
    \label{tab:dimacs}
\end{sidewaystable}

\begin{table}
    \centering
    \setlength{\tabcolsep}{0.1em}
    \scriptsize
    \begin{tabular}{ll@{\hskip 0.5em}r@{\hskip 0.5em}rr@{\hskip 0.5em}rr@{\hskip 0.5em}rr@{\hskip 0.5em}rr@{\hskip 0.5em}rr}
\toprule
 \multicolumn{2}{c@{\hskip 0.5em}}{Graph} &  \multicolumn{1}{c@{\hskip 0.5em}}{$\theta^+(G)$} & \multicolumn{2}{c@{\hskip 0.5em}}{$\lambda(G)$} & \multicolumn{2}{c@{\hskip 0.5em}}{$\mu(G, \mathcal{C})$} & \multicolumn{2}{c@{\hskip 0.5em}}{$\nu(G, \Gamma)$} & \multicolumn{2}{c@{\hskip 0.5em}}{$\nu(G, \theta)$} & \multicolumn{2}{c}{$\nu(G, \alpha)$} \\
Name & $d$(\%) & Time & \# iter & Time & \# iter & Time & \# iter & Time & \# iter & Time & \# iter & Time \\
\midrule
\texttt{brock200\_1}  & 25 & 3.72 & 1 & 7.22 & 1 & 6.86 & 0 & 3.79 & 0 & 3.78 & 0 & 3.79 \\
\texttt{brock200\_2}  & 50 & 1.63 & 0 & 2.51 & 1 & 4.02 & 0 & 1.73 & 1 & 3.79 & \replace{1}{2} & 14.22 \\
\texttt{brock200\_3}  & 39 & 2.07 & 1 & 4.27 & 1 & 4.44 & 0 & 2.15 & 1 & 3.92 & 1 & 4.47 \\
\texttt{brock200\_4}  & 34 & 2.66 & 0 & 3.12 & 1 & 4.79 & 0 & 2.75 & 0 & 2.75 & 1 & 4.35 \\
\midrule
\texttt{brock400\_1}  & 25 & 11.59 & 0 & 15.10 & 1 & 20.59 & 0 & 12.00 & 0 & 12.07 & 0 & 11.99 \\
\texttt{brock400\_2}  & 25 & 12.35 & 0 & 15.47 & 1 & 21.18 & 0 & 12.75 & 0 & 12.77 & 0 & 12.76 \\
\texttt{brock400\_3}  & 25 & 12.17 & 0 & 14.87 & 1 & 21.40 & 0 & 12.59 & 0 & 12.58 & 0 & 12.61 \\
\texttt{brock400\_4}  & 25 & 12.05 & 0 & 14.84 & 1 & 20.76 & 0 & 12.45 & 0 & 12.44 & 0 & 12.44 \\
\midrule
\texttt{brock800\_1}  & 35 & 20.60 & 0 & 51.54 & 0 & 27.57 & 0 & 24.54 & 0 & 24.60 & 2 & 77.95 \\
\texttt{brock800\_2}  & 35 & 20.26 & 0 & 53.58 & 0 & 27.05 & 0 & 24.25 & 0 & 24.22 & 2 & 86.15 \\
\texttt{brock800\_3}  & 35 & 20.55 & 0 & 52.85 & 0 & 28.19 & 0 & 25.02 & 0 & 24.63 & 2 & 79.55 \\
\texttt{brock800\_4}  & 35 & 19.86 & 0 & 50.82 & 0 & 27.46 & 0 & 23.98 & 0 & 24.00 & 2 & 76.29 \\
\midrule
\texttt{C125-9}       & 10 & 3.46  & 2 & 15.68 & \replace{3}{4} & 32.44 & 1 & 7.10  & 1 & 7.44  & 1 & 7.55 \\
\texttt{C250-9}       & 10 & 10.45 & \replace{1}{2} & 32.52 & 2 & 42.68 & 0 & 10.52 & 0 & 10.51 & 0 & 10.52 \\
\texttt{C500-9}       & 10 & 25.12 & 1 & 57.41 & 1 & 53.55 & 0 & 25.52 & 0 & 25.50 & 0 & 25.50 \\
\midrule
\texttt{DSJC125.1}    &  9 & 4.38  & 2 & 19.28 & \replace{3}{4} & 34.02 & 1 & 8.38  & 1 & 10.27 & 1 & 9.02 \\
\texttt{DSJC125.5}    & 50 & 1.58  & 1 & 3.17  & 1 & 3.72  & 0 & 1.61  & 1 & 4.35  & 1 & 4.64 \\
\texttt{DSJC125.9}    & 90 & 1.65  & 0 & 1.98  & 0 & 1.72  & 0 & 1.70  & 0 & 1.70  & 0 & 1.70 \\
\texttt{DSJC500-5}    & 50 & 6.35  & 0 & 20.49 & 0 & 8.85  & 0 & 7.71  & 0 & 7.83  & 4 & 126.12 \\
\midrule
\texttt{MANN\_a9 }    &  7 & 0.79  & 1 & 2.15  & \replace{1}{2} & 4.38    & 0 & 0.80  & 1 & 2.20  & 1 & 1.82 \\
\texttt{MANN\_a27 }   &  1 & 9.93  & 1 & 73.41 & 3 & 1992.72 & 0 & 10.00 & 1 & 57.22 & 1 & 60.32 \\
\texttt{johnson32-2-4}& 12 & 4.16  & 0 & 6.92  & 0 & 4.28    & 0 & 4.58  & 0 & 4.58  & 0 & 4.60 \\
\texttt{keller4}      & 35 & 3.97  & 1 & 7.70  & 1 & 13.20   & 0 & 4.03  & 1 & 15.51 & 1 & 15.12 \\
\midrule
\texttt{p\_hat300-1}  & 76 & 9.07  & 0 & 12.62  & 1 & 22.93  & 0 & 9.49  & 2 & 134.81 & 3 & 261.52 \\
\texttt{p\_hat300-2}  & 51 & 80.54 & \replace{1}{2} & 220.85 & 2 & 221.98 & 0 & 80.85 & 1 & 300.55 & 1 & 279.48 \\
\texttt{p\_hat300-3}  & 26 & 17.66 & 1 & 32.39  & 2 & 44.90  & 0 & 17.84 & 1 & 39.16  & 1 & 39.27 \\
\midrule
\texttt{p\_hat500-1}  & 75 & 19.08 & 0 & 40.72  & 1 & 42.03  & 0 & 21.13  & 2 & 164.90 & 5 & 858.88 \\
\texttt{p\_hat500-2}  & 50 & 201.27& 2 & 661.29 & 2 & 652.87 & 0 & 202.66 & 1 & 580.22 & 1 & 893.64 \\
\texttt{p\_hat500-3}  & 25 & 35.75 & 1 & 88.95  & 2 & 116.21 & 0 & 36.50  & 1 & 79.32  & 1 & 90.28 \\
\midrule
\texttt{p\_hat700-1}  & 75 & 35.53 & 0 & 82.10   & 1 & 79.12   & 0 & 41.34  & 2 & 260.41 & 7 & 6580.82 \\
\texttt{p\_hat700-2}  & 50 & 425.90& 2 & 1270.35 & 2 & 1409.98 & 0 & 429.64 & 1 & 920.86 & 1 & 1251.71 \\
\texttt{p\_hat700-3}  & 25 & 81.42 & 1 & 204.92  & 2 & 278.70  & 0 & 83.93  & 1 & 189.04 & 1 & 204.45 \\
\midrule
\texttt{sanr200\_0.7} & 30 & 3.05  & 1 & 5.39  & 1 & 5.54  & 0 & 3.11 & 0 & 3.11  & 0 & 3.11 \\
\texttt{sanr200\_0.9} & 10 & 7.77  & 2 & 29.39 & 3 & 43.35 & 0 & 7.81 & 1 & 15.36 & 1 & 14.86 \\
\texttt{sanr400\_0.5} & 50 & 4.40  & 0 & 11.74 & 1 & 8.96  & 0 & 5.09 & 0 & 5.11  & 3 & 77.60 \\
\texttt{sanr400\_0.7} & 30 & 7.69  & 0 & 11.58 & 1 & 12.80 & 0 & 8.16 & 0 & 8.17  & 0 & 8.26 \\
\bottomrule
\end{tabular}

    \caption{CPU time and cutting plane's iterations on DIMACS graphs.}
    \label{tab:dimacs_times}
\end{table}

\begin{table}
    \centering
        \setlength{\tabcolsep}{0.1em}
    \begin{tabular}{l@{\hskip 2.5em}r@{\hskip 2.5em}rr@{\hskip 2.5em}r@{\hskip 2.5em}rr@{\hskip 0.3in}rr}
\toprule
 \multicolumn{9}{c}{Number of violated inequalities added in Algorithm~\ref{alg:separation}} \\
\toprule
 \multicolumn{1}{c@{\hskip 2.5em}}{Graph}& \multicolumn{1}{c@{\hskip 2.5em}}{$\lambda(G)$} & \multicolumn{2}{c@{\hskip 2.5em}}{$\mu(G, \mathcal{C})$}  & \multicolumn{1}{c@{\hskip 2.5em}}{$\nu(G, \Gamma)$} & \multicolumn{2}{c@{\hskip 0.3in}}{$\nu(G, \theta)$} & \multicolumn{2}{c}{$\nu(G, \alpha)$} \\
\multicolumn{1}{c@{\hskip 0.3in}}{Name} & \eqref{eq:LS1} & \eqref{clq2}& \eqref{clq4} & \eqref{nod3}& \eqref{nod2}& \eqref{nod3}& \eqref{nod2}& \eqref{nod3}\\
\midrule
\texttt{brock200\_1   } &    10 &   471 &      0 &       0 &      0 &      0 &      0 &     0   \\
\texttt{brock200\_2   } &     0 &    81 &      0 &       0 &      5 &      4 &    583 &    76   \\
\texttt{brock200\_3   } &     1 &   152 &      0 &       0 &      1 &      1 &      2 &     5   \\
\texttt{brock200\_4   } &     0 &   139 &      0 &       0 &      0 &      0 &      1 &     0   \\
\texttt{brock400\_1   } &     0 &    34 &      0 &       0 &      0 &      0 &      0 &     0   \\
\texttt{brock400\_2   } &     0 &    38 &      0 &       0 &      0 &      0 &      0 &     0   \\
\texttt{brock400\_3   } &     0 &    49 &      0 &       0 &      0 &      0 &      0 &     0   \\
\texttt{brock400\_4   } &     0 &    33 &      0 &       0 &      0 &      0 &      0 &     0   \\
\texttt{brock800\_1   } &     0 &     0 &      0 &       0 &      0 &      0 &   1845 &   155   \\
\texttt{brock800\_2   } &     0 &     0 &      0 &       0 &      0 &      0 &   1780 &   220   \\
\texttt{brock800\_3   } &     0 &     0 &      0 &       0 &      0 &      0 &   1825 &   175   \\
\texttt{brock800\_4   } &     0 &     0 &      0 &       0 &      0 &      0 &   1840 &   160   \\
\texttt{C125-9        } &  1413 &  2248 &      3 &       1 &      0 &     35 &      0 &    14    \\
\texttt{C250-9        } &   565 &  1988 &      0 &       0 &      0 &      0 &      0 &     0    \\
\texttt{C500-9        } &     1 &   633 &      0 &       0 &      0 &      0 &      0 &     0    \\
\texttt{DSJC125.1     } &  1435 &  2249 &     11 &       8 &      0 &     36 &      0 &    36    \\
\texttt{DSJC125.5     } &     8 &   210 &      0 &       0 &     62 &     17 &    117 &    24    \\
\texttt{DSJC125.9     } &     0 &     0 &      0 &       0 &      0 &      0 &      0 &     0    \\
\texttt{DSJC500-5     } &     0 &     0 &      0 &       0 &      0 &      0 &   3772 &   228    \\
\texttt{MANN\_a9      } &   720 &   900 &     36 &       0 &      0 &    504 &      0 &   504   \\
\texttt{MANN\_a27     } &  1000 &  2000 &   1000 &       0 &      0 &   1000 &      0 &  1000   \\
\texttt{johnson32-2-4 } &     0 &     0 &      0 &       0 &      0 &      0 &      0 &     0    \\
\texttt{keller4       } &    48 &   296 &      0 &       0 &     96 &     48 &     96 &    48    \\
\texttt{p\_hat300-1   } &     0 &    48 &      0 &       0 &   1911 &     89 &   2923 &    77   \\
\texttt{p\_hat300-2   } &   988 &  1379 &      0 &       0 &    121 &    346 &    329 &   398   \\
\texttt{p\_hat300-3   } &   761 &  1878 &      0 &       0 &      0 &     60 &      0 &    75   \\
\texttt{p\_hat500-1   } &     0 &    13 &      0 &       0 &   1982 &     18 &   4959 &    41   \\
\texttt{p\_hat500-2   } &  1042 &  1363 &      0 &       0 &     30 &    288 &    576 &   424   \\
\texttt{p\_hat500-3   } &   548 &  1625 &      0 &       0 &      0 &     28 &      0 &    72   \\
\texttt{p\_hat700-1   } &     0 &     2 &      0 &       0 &   1960 &     40 &   6962 &    38   \\
\texttt{p\_hat700-2   } &  1107 &  1549 &      0 &       0 &     16 &    286 &    718 &   282   \\
\texttt{p\_hat700-3   } &   520 &  1414 &      0 &       0 &      0 &     47 &      1 &   187   \\
\texttt{sanr200\_0.7  } &     3 &   286 &      0 &       0 &      0 &      0 &      0 &     0   \\
\texttt{sanr200\_0.9  } &  1060 &  2216 &      0 &       0 &      0 &      1 &      0 &     1   \\
\texttt{sanr400\_0.5  } &     0 &     3 &      0 &       0 &      0 &      0 &   2806 &   194   \\
\texttt{sanr400\_0.7  } &     0 &    29 &      0 &       0 &      0 &      0 &      0 &     0   \\
\bottomrule
\end{tabular}

    \caption{Number of added cutting planes on DIMACS instances.}
    \label{tab:dimacs_cuts}
\end{table}

\add{\subsection{Sensitivity to the choice of $\mathcal{C}$}
To enhance the results of our computational analysis, we investigate whether refining the edge clique cover $\mathcal{C}$ obtained from Algorithm~\ref{alg:greedy_clq_cover} can lead to a significant improvement in the bound $\mu(G, \mathcal{C})$, while also considering the associated computational overhead.} 

\add{To this end, we consider the edge clique cover $\mathcal{C_{\rm cut}} \in \mathcal{Q}$  proposed by \citet{LeRoSm2020}. The collection $\mathcal{C_{\rm cut}}$ is obtained by adding to $\mathcal{C}$ all cliques identified through a cutting plane algorithm that employs the separation heuristic of~\citet{MaRoSm2019}. Since $\mathcal{C} \subseteq \mathcal{C_{\rm cut}}$,  $\ON{QSTAB}(G, \mathcal{C_{\rm cut}})$ provides a tighter upper bound on $\alpha(G)$ than that of $\qstabc{G}$.  \citet{MaRoSm2019} and~\citet{LeRoSm2020} also show that this upper bound is very close to the bound provided by $\qstab{G}$. However, this improvement may require a large number of clique inequalities.}

\add{Tables~\ref{tab:dimacs_sensitivity} and~\ref{tab:dimacs_sensitivity_times} report the results of this sensitivity analysis on a subset of the previously tested DIMACS instances. We selected graphs of varying densities for which we expect that $\mathcal{C_{\rm cut}}$ has the potential to 
overturn the conclusion drawn from the previous experiment. As anticipated, the percentage gaps in Table~\ref{tab:dimacs_sensitivity} show an improvement of $\mu(G, \mathcal{C_{\rm cut}})$ over 
$\mu(G, \mathcal{C})$. This gain, however, generally comes at a substantial additional cost in terms of cutting-plane iterations and CPU time (up to a factor of 2), as reported in Table~\ref{tab:dimacs_sensitivity_times}.
Interestingly, $\mathcal{C_{\rm cut}}$ has a negligible effect on the bounds for dense graphs, where $\nu(G, \alpha)$ remains superior. Overall, this experiment suggests that the greedy edge clique cover provides a good trade-off between bound quality and computational time.}

\add{A notable exception is \texttt{keller4}, where $\mu(G, \mathcal{C}_{\rm cut})$ attains a tight upper bound of $12.64$ on $\alpha(G)=11$, computed in a reasonable time. To the best of the authors’ knowledge, this represents the strongest bound reported in the literature; see, e.g.,~\citet{Lo2015, GiRoSm2013} and references therein.}

\section{Conclusions}
\label{sec:conclusion}
We have presented a way to improve the classical Lov\'asz $\theta$ bound for the stability number of a graph by applying the Lov\'asz-Schrijver \N{} lift-and-project operator to tailored LP relaxations. Unlike previous approaches, some of the latter may not be constructed in polynomial time. In the case of clique inequalities, we have to handle exponentially many inequalities (with an associated NP-hard separation problem), while building strong nodal inequalities requires solving the SSP on specific subgraphs. We show that these expedients, if properly handled, are indeed helpful in letting the resulting SDP relaxations improve the $\theta$ bound without excessive extra effort. This study reveals that the \lov{} and \sch{} \add{\N{}} operator\remove{s} can be of practical interest besides being a powerful theoretical tools. 
\add{Therefore, its application to other combinatorial optimization problems that admit strong SDP relaxations, such as the graph coloring problem, is a natural research line.}
\add{From a theoretical standpoint, a few questions remain open. 
First, from Theorem~\ref{th:theta} we have $\Nnod{G}{\Gamma} \subseteq \ON{TH}^+(G)$, whereas 
the numerical experiments suggest $\Nnod{G}{\Gamma} = \ON{TH}^+(G)$. However, it still remains an open question.
Secondly, the experiments showed a strong predominance of violated constraints lifted by multiplying by $x_i \ge 0$ on 
those multiplied by $1 - x_i \ge 0$. This deserves a further theoretical investigation.}
Finally, an interesting research direction consists of finding classes of graphs where $\Nnod{G}{\theta} = \ON{STAB}(G)$, as for them the SSP would be polynomially solvable.\\

\if\journal1
\begin{acknowledgements}
\else
\paragraph{Acknowledgements}
\fi
    We wish to thank Adam Letchford and Stefano Gualandi for helpful discussions about the preliminary results of this study.
\if\journal1
\end{acknowledgements}
\fi

\section*{Conflict of interest}
The authors declare that they have no conflict of interest.

\begin{table}
    \centering
        \setlength{\tabcolsep}{0.4em}

    \begin{tabular}{lrr@{\hskip 0.2in}rr@{\hskip 0.2in}rr@{\hskip 0.2in}rr}
\toprule
 \multicolumn{1}{c@{\hskip 0.2in}}{Graph} & 
 \multicolumn{2}{c@{\hskip 0.2in}}{$\theta^+(G)$} & 
 \multicolumn{2}{c@{\hskip 0.2in}}{$\mu(G, \mathcal{C_{\rm cut}})$} & 
 \multicolumn{2}{c@{\hskip 0.2in}}{$\mu(G, \mathcal{C})$} & 
 \multicolumn{2}{c}{$\nu(G, \alpha)$}  \\
 Name  & UB & \% gap & UB & \% gap & UB & \% gap & UB & \% gap\\
\midrule
\texttt{brock200\_1 }    & 27.20 & 29.508 & \textbf{26.71} & \textbf{27.183} & 27.12 & 29.119   & 27.20 & 29.508 \\
\texttt{brock200\_2 }    & 14.13 & 17.758 & \textbf{13.98} & \textbf{16.510} & 14.12 & 17.674   & 14.02 & 16.795 \\
\texttt{brock200\_3 }    & 18.67 & 24.479 & \textbf{18.30} & \textbf{22.013} & 18.65 & 24.356   & 18.67 & 24.467 \\
\texttt{brock200\_4 }    & 21.12 & 24.242 & \textbf{20.78} & \textbf{22.246} & 21.10 & 24.101   & 21.12 & 24.242 \\
\midrule
\texttt{C125-9      }    & 37.55 & 10.431 & \textbf{35.87} & \textbf{5.487} & 36.23 & 6.566          & 37.55 & 10.430 \\
\texttt{DSJC125.1   }    & 38.04 & 11.896 & \textbf{36.42} & \textbf{7.125} & 36.79 & 8.207       & 38.04 & 11.881 \\
\texttt{DSJC125.5   }    & 11.40 & 14.021 & \textbf{11.04} & \textbf{10.359} & 11.36 & 13.598     & 11.35 & 13.531 \\
\texttt{keller4     }    & 13.47 & 22.417 & \textbf{12.64} & \textbf{14.873} & 13.39 & 21.709       & 13.45 & 22.236 \\
\midrule
\texttt{p\_hat300-1 }    & 10.02 & 25.253 & 9.81 & 22.622 & 10.02 & 25.194                      & \textbf{8.58} & \textbf{7.288} \\
\texttt{p\_hat300-2 }    & 26.71 & 6.855 & \textbf{26.00} & \textbf{4.015} & 26.49 & 5.950      & 26.69 & 6.768 \\
\texttt{p\_hat300-3 }    & 40.70 & 13.056 & \textbf{39.74} & \textbf{10.384} & 40.44 & 12.321   & 40.70 & 13.053 \\
\midrule
\texttt{p\_hat500-1 }    & 13.01 & 44.533 & 12.95 & 43.847 & 13.01 & 44.526                     & \textbf{10.68} & \textbf{18.721} \\
\texttt{p\_hat500-2 }    & 38.56 & 48.306 & \textbf{37.94} & \textbf{45.934} & 38.38 & 47.597   & 38.54 & 48.225 \\
\midrule
\texttt{p\_hat700-1 }    & 15.05 & 36.774 & 15.02 & 36.589 & 15.05 & 36.773                     & \textbf{11.30} & \textbf{2.709} \\
\texttt{p\_hat700-2 }    & 48.44 & 10.091 & \textbf{47.89} & \textbf{8.831} & 48.24 & 9.647     & 48.41 & 10.023 \\
\midrule
\texttt{sanr200\_0.7}    & 23.63 & 31.296 & \textbf{23.22} & \textbf{28.976} & 23.59 & 31.048  & 23.63 & 31.296 \\
\texttt{sanr200\_0.9}    & 48.90 & 16.439 & \textbf{47.49} & \textbf{13.061} & 48.03 & 14.356  & 48.90 & 16.440 \\
\texttt{sanr400\_0.5}    & 20.18 & 55.217 & 20.17 & 55.142 & 20.18 & 55.216                    & \textbf{19.10} & \textbf{46.951} \\
\bottomrule
\end{tabular}

    \caption{Semidefinite bounds with improved edge clique covers.}
    \label{tab:dimacs_sensitivity}
\end{table}

\begin{table}
    \centering
    \begin{tabular}{l@{\hskip 0.2in}rr@{\hskip 0.2in}rr@{\hskip 0.2in}rr@{\hskip 0.2in}rr@{\hskip 0.2in}rr}
\toprule
 \multicolumn{1}{c@{\hskip 0.2in}}{Graph} 
 & \multicolumn{2}{c@{\hskip 0.2in}}{$\mu(G, \mathcal{C_{\rm cut}})$} 
 & \multicolumn{2}{c@{\hskip 0.2in}}{$\mu(G, \mathcal{C})$} 
 & \multicolumn{2}{c}{$\nu(G, \alpha)$} \\
Name  & \# iter & Time & \# iter & Time & \# iter & Time  \\
\midrule
\texttt{brock200\_1}  & 2.0 & 10.60 & 1.0 & 6.31        & 0.0 & 3.34 \\
\texttt{brock200\_2}  & 2.0 & 7.04 & 1.0 & 3.87         & 2.0 & 13.88 \\
\texttt{brock200\_3}  & 2.0 & 7.03 & 1.0 & 4.18         & 1.0 & 4.26 \\
\texttt{brock200\_4}  & 2.0 & 7.29 & 1.0 & 4.56         & 1.0 & 4.54 \\
\midrule
\texttt{C125-9     }  & 4.0 & 32.91 & 4.0 & 33.90            & 1.0 & 9.95 \\
\texttt{DSJC125.1  }  & 4.0 & 32.81 & 4.0 & 33.64         & 1.0 & 9.22 \\
\texttt{DSJC125.5  }  & 2.0 & 5.53 & 1.0 & 2.97           & 1.0 & 3.64 \\
\texttt{keller4    }  & 4.0 & 29.15 & 1.0 & 11.57           & 1.0 & 15.01 \\
\midrule
\texttt{p\_hat300-1} & 2.0 & 36.83 & 1.0 & 21.78       & 3.0 & 262.60 \\
\texttt{p\_hat300-2}  & 3.0 & 311.08 & 2.0 & 168.35    & 1.0 & 232.25 \\
\texttt{p\_hat300-3}  & 4.0 & 75.37 & 2.0 & 43.34      & 1.0 & 36.93 \\
\midrule
\texttt{p\_hat500-1} & 1.0 & 37.54 & 1.0 & 40.53      & 5.0 & 854.09 \\
\texttt{p\_hat500-2}  & 3.0 & 841.79 & 2.0 & 588.09   & 1.0 & 927.59 \\
\midrule
\texttt{p\_hat700-1} & 1.0 & 76.46 & 1.0 & 78.54      & 7.0 & 6307.61 \\
\texttt{p\_hat700-2}  & 3.0 & 1405.20 & 2.0 & 1332.51 & 1.0 & 1208.12 \\
\midrule
\texttt{sanr200\_0.7} & 2.0 & 9.51 & 1.0 & 5.50        & 0.0 & 3.19 \\
\texttt{sanr200\_0.9}  & 3.0 & 43.39 & 3.0 & 43.01     & 1.0 & 20.16 \\
\texttt{sanr400\_0.5} & 1.0 & 8.78 & 1.0 & 8.70        & 3.0 & 76.85 \\
\bottomrule
\end{tabular}

    \caption{CPU time and cutting-plane iterations with improved edge clique covers.}
    \label{tab:dimacs_sensitivity_times}
\end{table}

\clearpage

\bibliographystyle{abbrvnat} 
\bibliography{ref}

\end{document}